\documentclass[a4paper,12pt]{article}
\usepackage[utf8]{inputenc}

\title{Staggered grids for multidimensional multiscale modelling}

\author{
J. Divahar
\thanks{School of Mathematical Sciences,
University of Adelaide, South Australia.}
\thanks{\protect\url{https://orcid.org/0000-0002-9506-8846}}
\and
A.~J. Roberts
\footnotemark[1]
\thanks{\protect\url{http://orcid.org/0000-0001-8930-1552},
\protect\url{mailto:profajroberts@protonmail.com}}
\and
Trent W. Mattner
\footnotemark[1]
\thanks{\protect\url{https://orcid.org/0000-0002-5313-5887}}
\and
J.~E. Bunder
\footnotemark[1]
\thanks{\protect\url{http://orcid.org/0000-0001-5355-2288}}
\and
Ioannis~G. Kevrekidis
\thanks{Departments of Chemical and Biomolecular Engineering \& Applied Mathematics and Statistics, Johns Hopkins University, Baltimore, Maryland, USA.
\protect\url{https://orcid.org/0000-0003-2220-3522}}
}

\usepackage{amsfonts,amsmath,amssymb,amsthm}

\IfFileExists{ajr.sty}{
  \usepackage[a5paper,margin=5.5mm,tmargin=13mm]{geometry}
  \renewcommand{\vec}[1]{\text{\boldmath$##1$}}
  \AtBeginDocument{
    \setNodeShape{tri}
    \setuptodonotes{inline,color=orange!25}
    \addbibresource{ajr,bib}
  }
  \usepackage{wrapfig}

}{ 
  \usepackage[small, euler-digits]{eulervm} 
  \usepackage{palatino}  
  \let\vec\mathbf
  \AtBeginDocument{
    \setNodeShape{circ}  
    \setuptodonotes{inline,color=orange!25}
  }
  
}

\usepackage[fit]{truncate}
\usepackage{fancyhdr}
\usepackage{microtype} 
\usepackage{needspace}

\usepackage[svgnames]{xcolor}

\usepackage{graphicx}
\graphicspath{{figs/inkscape/}}

\usepackage{subcaption}

\usepackage{hyperref}
\usepackage{hypcap}
\hypersetup{
  pdftex,
  backref=true,
  hyperindex=true,
  colorlinks=true,
  allcolors=teal,
  bookmarks=true,
  breaklinks=true
}

\usepackage[capitalise,nameinlink,noabbrev]{cleveref}
\crefname{equation}{}{} 
\crefname{enumi}{}{}
\crefname{enumii}{}{}
\crefname{enumiii}{}{}
\crefname{enumiv}{}{}

\usepackage{citeSettings}
\makeatletter
\let\L@TeX@startsection\@startsection
\renewcommand{\@startsection}[6]{\L@TeX@startsection{#1}{#2}{#3}{#4}{#5}
{\raggedright #6}}
\makeatother

\usepackage{defs}

\usepackage[textwidth=3.3cm,textsize=small]{todonotes}
\usepackage{marginnote}

\def\ZZ{\mathbb Z}
\def\ev{\vec{e}}

\begin{document}

\maketitle

\begin{abstract}
Numerical schemes for wave-like systems with small dissipation are often inaccurate and unstable due to truncation errors and numerical roundoff errors.
Hence, numerical simulations of wave-like systems lacking proper handling of these numerical issues often fail to represent the physical characteristics of wave phenomena.
This challenge gets even more intricate for multiscale modelling, especially in multiple dimensions.
%
When using the usual collocated grid, about two-thirds of the resolved wave modes are incorrect with significant dispersion.
But, numerical schemes on staggered grids (with alternating variable arrangement) are significantly less dispersive and preserve much of the wave characteristics.
Also, the group velocity of the energy propagation in the numerical waves on a staggered grid is in the correct direction, in contrast to the collocated grid.
For high accuracy and to preserve much of the wave characteristics, \emph{this article extends the concept of staggered grids in full-domain modelling to multidimensional multiscale modelling}.
Specifically, this article develops \(120\) multiscale staggered grids and demonstrates their stability, accuracy, and wave-preserving characteristic for equation-free multiscale modelling of weakly damped linear waves.
But most characteristics of the developed multiscale staggered grids must also hold in general for multiscale modelling of many complex spatio-temporal physical phenomena such as the general computational fluid dynamics.
\end{abstract}

\tableofcontents

\section{Introduction}

For wave-like systems with small or no dissipation, accurate numerical simulation is challenging over large spatial scales, especially for long simulation times.
Numerical schemes for wave-like systems with small dissipation are often inaccurate and unstable due to numerical dissipation and numerical dispersion caused by truncation and numerical roundoff errors (\cites[p.136]{Hinch2020_ThnkBfreYuCmpte_APrldeToCmpttnlFldDynmcs}[pp.70--73]{Zikanov2010_EsntlCmpttnlFldDynmcs}[pp.232--243]{Anderson1995_CmpttnlFldDynmcs}).
Hence, numerical simulations of wave-like systems lacking proper handling of these numerical issues often fail to represent the physical characteristics of wave phenomena.
\Cref{sec:cllctdAndStgrdGrdsFrWvelkeSystms} overviews using a staggered grid in space for accurate and robust computational simulation of wave-like systems.
It is well-known that staggered grids, such as depicted in \cref{fig:flDmnGrd_clctd_stgrd}(right), lead to higher accuracy compared to that of a same order scheme \text{on collocated grids.}

Here we develop multiscale staggered grids for wave-like systems in multiple space dimensions---such systems are even more challenging (e.g., see the recent review of modelling materials by \cite{Fish2021}).
We expect that for simulations of wave-like systems many of the good characteristics of the usual staggered grids also hold for the multiscale staggered grids.
The \emph{Equation-Free Patch scheme} \cite[e.g.]{Kevrekidis2009_EqtnFreMltscleCmpttnAlgrthmsAndAplctns} is a flexible, computationally efficient, multiscale modelling approach.
Such patch schemes are an effective and accurate alternative to cognate methods such as the Multiscale Finite Element Method \cite[e.g.]{Efendiev2004, Fish2021, Altmann2021}.
Equation-free multiscale patch schemes have been developed and applied successfully for dissipative systems \parencite{Roberts2005_HiOrdrAcrcyGapTthScme,Roberts2007_GnrlTthBndryCndnsEqnFreeMdlng,Bunder2017_GdCplingFrTheMltsclePtchSchmeOnSystmsWthMcrscleHtrgnty,Maclean2021_AtlbxOfEqtnFreFnctnsInMtlbOctveFrEfcntSystmLvlSmltn}.
\cite{Cao2013_MltSclMdlCplsWave, Cao2015_MltSclMdlCplsNonLrWave} extended the patch scheme to 1D wave-like systems using a staggered macroscale grid of patches in 1D space, where each patch itself contains a staggered microscale grid in 1D space.
Here, \cref{sec:stgrdPtchGrdsFrEqtnFreMltscleMdlngOfWvs} extends and explores the multiscale staggered grids to two-dimensional multiscale modelling (\cref{ssc:extndTheStgrdGrdToMltsclePtchSchme}).  
Extending to higher dimensional space should be straightforward.
For wave-like systems in 2D, a total of~\(167\,040\) such staggered patch grids are \emph{geometrically compatible} (defined in \cref{ssc:extndTheStgrdGrdToMltsclePtchSchme}).
We developed and analysed all \(167\,040\) staggered patch grids, among which \(120\) of them are \emph{centred staggered patch grids} whose sub-patch micro-grids all have a centre node and include a centred-patch corresponding to each of the state variables; the remaining \(166\,920\) are \emph{non-centred patch grids}.
For first-order \pde{}s, the patch schemes over these staggered patch grids interpolate field values to the patch edges.
However, higher-order spatial derivatives (e.g., for diffusion) and some nonlinear terms require interpolation to additional layers of edge values (\cref{ssc:mltiLyrEdgeNdsFrHghrOrdrSptlDrvtvs}).

\Cref{sec:mltscleStgrdPtchGrdFrWklyDmpdLnrWvs} shows that all the \(60\)~centred multiscale staggered patch grids give stable and accurate equation-free multiscale patch schemes.
For representative physical and discretisation parameters, via eigenvalues and wave frequencies, \cref{sec:mltscleStgrdPtchGrdFrWklyDmpdLnrWvs} demonstrates the stability, accuracy, and wave-preserving characteristic of these centred multiscale staggered grids for weakly damped linear waves.
\Cref{ssc:onlyPtchGrdsWthAlSymtrcPtchsAreStble} studies the dependence of the patch scheme stability on the patch grid geometry, the \emph{geometry-stability study}, for all the \(167\,040\) compatible 2D staggered patch grids.
The geometry-stability study shows that among all the compatible~\(167\,040\) compatible patch grids, only \(1248\) patch grids (\(0.75\%\)) whose patches are all symmetric are stable.
\Cref{ssc:onlyCntrdPtchGrdsAreAcrte} shows that none of the non-centred patch grids that are stable (\cref{ssc:onlyPtchGrdsWthAlSymtrcPtchsAreStble}) is accurate.
Thus, among all the possible~\(167\,040\) compatible 2D staggered patch grids, only \(120\) centred patch grids (\(0.07\%\)) are both stable and accurate.

Given the new foundation established here for multidimensional multiscale staggered grids in the Equation-Free Patch scheme, our subsequent article on multiscale simulation of large-scale linear waves explores two families of patch schemes in detail for their stability, accuracy, and consistency. 
That forthcoming article also verifies the schemes' insensitivity to numerical roundoff errors.
Subsequent articles are planned to develop and explore staggered patch schemes for nonlinear wave \pde{}s, specifically for the viscous and turbulent shallow water flows.
These developments aim to work towards novel, powerful, and efficient computational simulation of the complex multiscale physics of floods, bores, and tsunamis \cite[e.g.]{LeVeque2015, Reungoat2018}.

\section{Collocated and staggered grids for wave-like systems}
\label{sec:cllctdAndStgrdGrdsFrWvelkeSystms}

\begin{figure}
  \caption{\label{fig:flDmnGrd_clctd_stgrd}%
    Schematic \emph{micro-grids}:  \emph{collocated grid} on left; \emph{staggered grid} on right where variables are stored only on staggered/alternating discrete points (\emph{nodes}~\(\hNode\,h\), \(\uNode\,u\), and~\(\vNode\,v\)).
    The collocated grid has \(3 \times 3\) (\(n=3\)) grid intervals in the {\cij green} grid, whereas the staggered grid has \(6 \times 6\) (\(n=6\)) grid intervals, yet both have \(3 \times 3 = 9\) micro-cells ({\cpq orange} square).
    The unfilled nodes (\(\hNodeB\), \(\uNodeB\), and~\(\vNodeB\)) indicate discrete \(n\)-periodic boundary values.
  }
  \centering
  \small
  \input{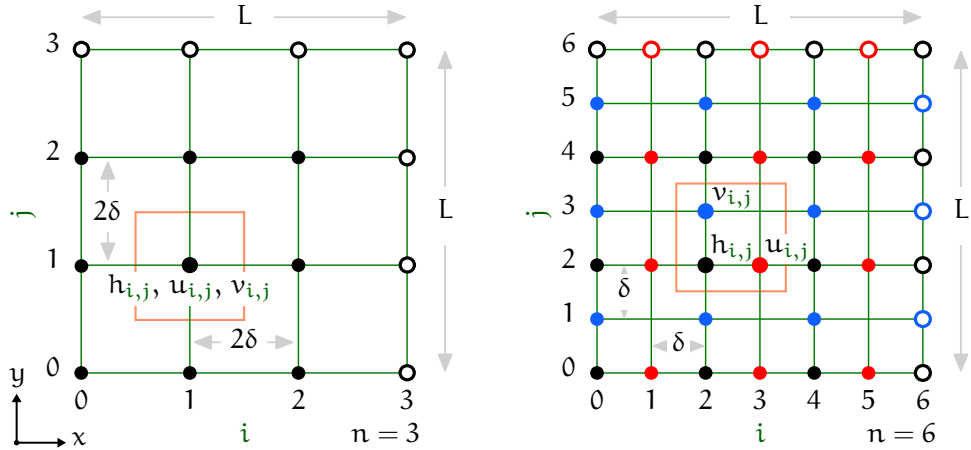}
\end{figure}

\Cref{fig:flDmnGrd_clctd_stgrd}(left) depicts an example collocated grid in 2D~space, and \cref{fig:flDmnGrd_clctd_stgrd}(right) a example 2D staggered grid.
In the classification of \textcite[Fig.~3, p.181]{Arakawa1977_CmpttnlDsgnOfTheBscDynmclPrcssOfTheUCLAGnrlCrcltnMdl}, these grids are A-grid and C-grid respectively.
The collocated grid in \cref{fig:flDmnGrd_clctd_stgrd}(left) stores all state variables at each and every discrete point, also called  node, at each~\((i,j)\) where {\cij green} grid lines intersect.
But the staggered grid in \cref{fig:flDmnGrd_clctd_stgrd}(right) intersperses the state variables at alternate discrete points.
Here, \(h, u\)~nodes are horizontally alternating and \(h, v\)~nodes are vertically alternating.
The staggered grid was first used in the Marker and Cell method of \textcite{Harlow1965_NmrclClcltonOfTmeDpndntVscsIncmprsbleFlwOfFldWthFreSrfce}.
Staggered grids were subsequently used in many contexts, including the well-known Semi Implicit Method for Pressure Linked Equations (\textsc{simple}) of \textcite{Patankar1972_AClcltnPrcdreFrHtMsAndMmntmTrnsFrInThreDmnsnlPrblcFlws}.

\needspace{3\baselineskip}
Staggered grids preserve much of the wave characteristics (\cites
  [Figs.~8 and 9]{Fornberg1999_SptlFnteDfrnceAprxmtnsFrWveTypeEqtns}%
  {Fornberg1990_HghOrdrFnteDfrncsAndThePsdspctrlMthdOnStgrdGrds}%
  ) and typically support higher accuracy simulations compared to simulations of the same order on collocated grids.
For example, even though a central difference scheme on a collocated grid gives second-order accuracy, the same second-order accurate central difference scheme on a staggered grid is significantly less dispersive (\cites%
  [p.46, \S3.2]{Lauritzen2011_NmrclTchnqsFrGlblAtmsphrcMdls}%
  [p.55, \S2.2.1]{Olafsson2021_UncrtntsInNmrclWthrPrdctn}%
  ).
Furthermore, the group velocity of the energy propagation in the numerical waves on a staggered grid is approximately in the correct direction, whereas on collocated grids the group velocity for large wavenumbers is often in the opposite direction (\cites%
[p.46, \S3.2]{Lauritzen2011_NmrclTchnqsFrGlblAtmsphrcMdls}%
[p.55, \S2.2.1]{Olafsson2021_UncrtntsInNmrclWthrPrdctn}%
).
When using the usual collocated grid, about two-thirds of the resolved wave modes are incorrect with significant dispersion.
For waves in multiple dimensions, such errors in dispersion and group velocity extend over a larger portion of wavenumber space, and so a staggered grid is significantly better than a corresponding scheme on a collocated grid.
Staggered grids, such as \cref{fig:flDmnGrd_clctd_stgrd}, are simple and robust for numerical simulation of wave-like systems.

We now describe second-order accurate, finite difference schemes for generic 2D wave-like \pde{}s on both collocated and staggered grids.
Consider such a system with little or no dissipation over the periodic spatial domain \( [0,L] \times [0, L]\).
In terms of state variables~\(h(x,y,t)\), \(u(x,y,t)\), and~\(v(x,y,t)\), such systems are modelled as non-dimensional wave-like \pde{}s
\begin{equation} \label{eqs:PDEs_gnrcWve}
  \doh h t =-\doh u x-\doh v y +\cdots, \quad
  \doh u t =-\doh h x +\cdots, \quad
  \doh v t =-\doh h y +\cdots,
\end{equation}
with the boundary conditions that the three non-dimensional fields~\(h\), \(u\), and~\(v\) are \(L\)-periodic in both the~\(x\) and~\(y\) directions.
In the generic wave-like \pde{}s~\cref{eqs:PDEs_gnrcWve}, \emph{the explicitly written terms (excluding ``\(\cdots\)'') model non-dispersive non-dissipative \emph{ideal} wave phenomena} (\cites[pp.136--137]{Dean1991_WtrWveMchncsFrEngnrsAndScntsts}[p.260]{Mehaute1976_AnIntro2HydrdynmcsWtrWvs}).
The ``\(\cdots\)'' indicates other linear and\slash or nonlinear terms that model additional physical phenomena that modify the wave, such as bed drag, viscous diffusion, turbulent mixing, and surface tension.
The dependent variables in \pde{}s~\cref{eqs:PDEs_gnrcWve} let us interpret the \pde{}s as a model of water waves with wave height~\(h\) and horizontal velocities~\(u,v\), but the importance of~\cref{eqs:PDEs_gnrcWve} is that it is a generic model of many 2D wave phenomena.

\Cref{fig:flDmnGrd_clctd_stgrd}(left) depicts the usual \emph{collocated grid}, whereas \cref{fig:flDmnGrd_clctd_stgrd}(right) depicts a \emph{staggered grid} for typical  simulations of generic wave-like systems modelled by \pde{}s~\cref{eqs:PDEs_gnrcWve}.
The grids in \cref{fig:flDmnGrd_clctd_stgrd} are termed  \emph{micro-grids} in order to distinguish them from the macro-grids created for our multiscale modelling (\cref{sec:stgrdPtchGrdsFrEqtnFreMltscleMdlngOfWvs}).
\Cref{fig:flDmnGrd_clctd_stgrd} shows the collocated grid with \(3 \times 3\) (\(n=3\)) grid intervals, and the staggered grid with \(6 \times 6\) (\(n=6\)) grid intervals, yet both these grids contain \(3 \times 3 = 9\) \emph{micro-cells} (one micro-cell is depicted by an {\cpq orange} square).
In both collocated and staggered grids, each micro-cell contains all three state variables~\(\hVarMu\), \(\uVarMu\), and~\(\vVarMu\).
In general, for any positive even number~\(n\), a collocated grid with \(n/2 \times n/2\) grid intervals and a staggered grid with \(n \times n\) grid intervals, have the same number of micro-cells, namely \(n/2 \times n/2 = n^2/4\), and the same number of state variables, namely~\(3n^2/4\).
The unfilled nodes (\(\hNodeB\), \(\uNodeB\), and~\(\vNodeB\)) on the boundaries of the two grids in \cref{fig:flDmnGrd_clctd_stgrd} schematically indicate discrete \(n\)-periodic boundary values.

Consider the collocated grid of \cref{fig:flDmnGrd_clctd_stgrd}(left) with \(n\)~grid intervals, each \(2\delta\)-wide in both~\(x\) and~\(y\) directions (e.g., \(n=3\) for the specifiic collocated grid of \cref{fig:flDmnGrd_clctd_stgrd}(left)).
Approximating the spatial derivatives in the wave-like \pde{}s~\cref{eqs:PDEs_gnrcWve} on the nodes of the collocated grid (e.g., filled markers of \cref{fig:flDmnGrd_clctd_stgrd}(left)) by central finite differences gives the \emph{collocated microscale model} (for \({\cij i,j} \in \{0, 1, 2, \ldots, n-1\}\))
\begin{subequations} \label{eqs:clctdFDEs_gnrcWve}
\begin{align}
  \hNode\; \frac{d h_{{\cij i,j}}}{d t} & =
     -\frac{u_{\cij i+1,j} - u_{\cij i-1,j}}{4 \delta}
     - \frac{v_{\cij i,j+1} - v_{\cij i,j-1}}{4 \delta}
     + \cdots
     \,,\\
  \hNode\; \frac{d u_{{\cij i,j}}}{d t} & =
    - \frac{h_{\cij i+1,j} - h_{\cij i-1,j}}{4 \delta}
    + \cdots
     \,,\\
  \hNode\; \frac{d v_{{\cij i,j}}}{d t} & =
    - \frac{h_{\cij i,j+1} - h_{\cij i,j-1}}{4 \delta}
    + \cdots
    \,.
\end{align}
\end{subequations}
Similarly, consider the staggered grid in \cref{fig:flDmnGrd_clctd_stgrd}(right) with \(n\)~grid intervals each \(\delta\)-wide in both~\(x\) and~\(y\) directions (e.g., \(n=6\) for the specific staggered grid of \cref{fig:flDmnGrd_clctd_stgrd}(right)).
Approximating the spatial derivatives in the generic wave-like \pde{}s~\cref{eqs:PDEs_gnrcWve} on the nodes (e.g., filled markers of \cref{fig:flDmnGrd_clctd_stgrd}(right)) by central finite differences gives the \emph{staggered microscale model} (for \({\cij i,j} \in \{0, 1, 2, \ldots, n-1\}\))
\begin{subequations} \label{eqs:stgrdFDEs_gnrcWve}
\begin{align}
  \hNode\; \frac{d h_{{\cij i,j}}}{d t} & =
     -\frac{u_{\cij i+1,j} - u_{\cij i-1,j}}{2 \delta}
     - \frac{v_{\cij i,j+1} - v_{\cij i,j-1}}{2 \delta}
     + \cdots
     \quad i,j\text{ even}
     \,,\\
  \uNode\; \frac{d u_{{\cij i,j}}}{d t} & =
    - \frac{h_{\cij i+1,j} - h_{\cij i-1,j}}{2 \delta}
    + \cdots
    \quad i\text{ odd,}\quad j\text{ even}
    \,,\\
  \vNode\; \frac{d v_{{\cij i,j}}}{d t} & =
    - \frac{h_{\cij i,j+1} - h_{\cij i,j-1}}{2 \delta}
    + \cdots
        \quad i\text{ even,}\quad j\text{ odd}
    \,.
\end{align}
\end{subequations}
The ``\(\cdots\)'' in the collocated and staggered grid models~\cref{eqs:clctdFDEs_gnrcWve,eqs:stgrdFDEs_gnrcWve} indicate discrete approximations of other linear and/or nonlinear terms corresponding to the ``\(\cdots\)'' in the generic wave-like \pde{}s~\cref{eqs:PDEs_gnrcWve}.
The coloured bullets~\(\hNode, \uNode, \vNode\) in the  microscale models~\cref{eqs:clctdFDEs_gnrcWve,eqs:stgrdFDEs_gnrcWve} indicate the~\(h\), \(u\), and~\(v\) nodes, respectively.

Corresponding to the periodic boundary conditions of the generic wave-like \pde{}s~\cref{eqs:PDEs_gnrcWve}, in the collocated and staggered grid models~\cref{eqs:clctdFDEs_gnrcWve,eqs:stgrdFDEs_gnrcWve} the three fields~\(h\), \(u\), \(v\) are \(n\)-periodic in both~\(\cij i\) and~\(\cij j\), where \(n = L/(2\delta)\) for the collocated grid and \(n = L/\delta\) for the staggered grid.

The \pde{}s~\cref{eqs:PDEs_gnrcWve} inspire the  microscale computational models~\cref{eqs:clctdFDEs_gnrcWve,eqs:stgrdFDEs_gnrcWve}.
%
But, the aim of our multiscale patch schemes is to accurately compute the solutions of the microscale computational model, not the solutions of the \pde{}s.
\label{pge:pdePrphrl}
Thus, how well the  microscale discretisations~\cref{eqs:clctdFDEs_gnrcWve,eqs:stgrdFDEs_gnrcWve} predict the solutions of the \pde{}s~\cref{eqs:PDEs_gnrcWve} is a peripheral issue for this article.

Arranging the state variables of the microsccale model~\cref{eqs:clctdFDEs_gnrcWve} defined on the \emph{nodes} of the collocated grid in \cref{fig:flDmnGrd_clctd_stgrd}(left), namely~\(\hNode\,\hVarMu\), \(\hNode\,\uVarMu\), and \(\hNode\,\vVarMu\), into a vector gives the \emph{state vector of the collocated microscale model}~\(\vec{x}\).
Similarly, arranging the state variables of the microscale model~\cref{eqs:stgrdFDEs_gnrcWve} defined on the \emph{nodes} of the staggered grid in \cref{fig:flDmnGrd_clctd_stgrd}(right),  namely~\(\hNode\,\hVarMu\), \(\uNode\,\uVarMu\), and \(\vNode\,\vVarMu\), into a vector gives the \emph{state vector of the staggered microscale model}~\(\vec{x}\).
In terms of the state vector~\(\vec{x}\), both the collocated and staggered microscale models~\cref{eqs:clctdFDEs_gnrcWve,eqs:stgrdFDEs_gnrcWve} are represented as a dynamical system by the \ode{}s 
\begin{equation} \label{eqs:dynSys_mdN_gnrcWve}
  \dby{\vec{x}}{t} = \vec{f}( \vec{x} ),
\end{equation}
where the microscale model~\(\vec{f}\) is a vector function encapsulating the spatially discrete model, such as the equations~\cref{eqs:clctdFDEs_gnrcWve} or \cref{eqs:stgrdFDEs_gnrcWve}, for the generic wave-like \pde{}s~\cref{eqs:PDEs_gnrcWve}.
A full-domain microscale simulation is performed by numerical time-integration of the \ode{}s~\cref{eqs:dynSys_mdN_gnrcWve}---the spatially discrete microscale model~\cref{eqs:clctdFDEs_gnrcWve} or \cref{eqs:stgrdFDEs_gnrcWve}---on the nodes of the  collocated\slash staggered  micro-grid (filled markers in \cref{fig:flDmnGrd_clctd_stgrd}).

\section{Staggered patch grids for equation-free multiscale modelling of waves}
\label{sec:stgrdPtchGrdsFrEqtnFreMltscleMdlngOfWvs}

This section extends the concept of 2D staggered grids (\cref{sec:cllctdAndStgrdGrdsFrWvelkeSystms}) to multiscale modelling by developing \emph{2D staggered patch grids} for equation-free multiscale modelling \cite[e.g.]{Kevrekidis2009_EqtnFreMltscleCmpttnAlgrthmsAndAplctns, Cao2015_MltSclMdlCplsNonLrWave}.
A further extension to 3D should be straightforward.

In the equation-free multiscale framework, a patch scheme performs detailed microscale simulations within small, widely separated, \emph{patches} of space (e.g., the small {\cIJ violet} squares enclosing {\cij green} grid in \cref{fig:PtchGrd_n1t0_huvx}), and couples the patches (\emph{patch coupling}) via interpolation over the macroscale space between the patches \parencite{
  Kevrekidis2004_EqFreeCompAiddAnlysMltiSclSys,
  Kevrekidis2009_EqtnFreMltscleCmpttnAlgrthmsAndAplctns,
  Hyman2005_PtchDynmcsForMltsclePrblms}.
One can achieve arbitrarily high orders of macroscale consistency for patch schemes via appropriate high order interpolation for the patch coupling \parencite{
  Roberts2005_HiOrdrAcrcyGapTthScme,
  Roberts2007_GnrlTthBndryCndnsEqnFreeMdlng,
  Bunder2020a}.
For development and analysis of the staggered patch grids, throughout this article we use \emph{spectral interpolation} \parencite[\S2.2.3]{Bunder2020_LrgeScleSmltnOfShlwWtrWvsVaCmpttnOnlyOnSmlStgrdPtchs}, via Fourier transforms, to couple patches.
%

\Cref{ssc:extndTheStgrdGrdToMltsclePtchSchme} explains extending 2D staggered grids over a full-domain to multiscale modelling, using a 2D staggered macroscale grid of patches and a 2D staggered \emph{micro-grid} within each patch.
The discussion in \cref{ssc:extndTheStgrdGrdToMltsclePtchSchme} is for \pde{}s with only first order derivatives.
\Cref{ssc:mltiLyrEdgeNdsFrHghrOrdrSptlDrvtvs} explains extending the 2D staggered grid to multiscale modelling for linear\slash{}nonlinear \pde{}s with higher-order and mixed derivatives.
\cref{sec:mltscleStgrdPtchGrdFrWklyDmpdLnrWvs} then determines which of the many possible staggered patch grids give stable and accurate patch schemes for linear waves~\eqref{eqs:PDEs_gnrcWve}.

\subsection{Patch coupling connects the scales}
\label{ssc:ptchCplngCnctsTheScls}

In an equation-free multiscale patch scheme, patch coupling gives the patch edge values (unfilled circles on the patch edges in \cref{fig:PtchGrd_n1t0_huvx}), via two steps:
\begin{enumerate}
  \item  for each patch compute a representative \emph{macroscale patch value} (also called aggregate value, amplitude or order parameter) from the patch's respective microscale interior values;
  \item  compute the \emph{microscale} edge values of each patch by interpolating from the macroscale patch values of neighbouring patches across the relatively large inter-patch distances.
\end{enumerate}
In this way, patch coupling provides a two-way connection between the microscale and macroscale.

Patch coupling using polynomial or spectral interpolation leads to two broad families of patch schemes.
To focus on the issue of the geometric design of the multiscale patch scheme, throughout this article we use the highly accurate spectral interpolation to provide the edge values \cite[\S2.2.3]{Bunder2020_LrgeScleSmltnOfShlwWtrWvsVaCmpttnOnlyOnSmlStgrdPtchs}.
In summary, at each time one first computes the 2D Fourier transform of the patch's aggregate values.  
Second, for each edge position, simultaneously over all patches, one uses the Fourier shifting property, and an inverse 2D Fourier transform, to compute all the required edge values.
This spectral patch coupling was found to be highly accurate in \emph{one} possible design of a 2D staggered patch scheme \cite[\S2.2.4]{Bunder2020_LrgeScleSmltnOfShlwWtrWvsVaCmpttnOnlyOnSmlStgrdPtchs}.

In a patch scheme, at each time step, the edge values (by the patch coupling) and the interior node values (from the previous time step) are known.
Given the edge values and the interior values at the current time step, an \ode\ integrator such as \texttt{BS3} of \texttt{DifferentialEquations.jl} \parencite{Rackauckas2017_DfrntlEqtns.jl_APrfrmntAndFtreRchEcsystmFrSlvngDfrntlEqtnsInJla}, computes the interior values at the next time step by computing time derivatives from the governing \ode{}s, such as the \ode{}s~\cref{eqs:stgrdFDEs_gnrcWve}.

Patch schemes essentially provide a reduced order multiscale model of the given corresponding microscale model.
The macroscale (aggregate) values of patches would then usually be state variables in a slow manifold of the patch scheme (\cites
  [\S5.3, p.302]{Roberts2003_LwDmnsnlMdlngOfDynmclSystmsApldToSmeDsptveFldMchncs}
  [p.1349]{Kevrekidis2004_EqFreeCompAiddAnlysMltiSclSys}
  []{Zagaris2009_AnlyssOfTheAcrcyAndCnvrgnceOfEqtnFrePrjctnToASlwMnfld}
  []{Roberts1988_TheAplctnOfCntreMnfldThryToTheEvltnOfSystmWhchVrySlwlyInspce}
  []{Foias1988_OnTheCmpttnOfInrtlMnflds}
  []{Lorenz1986_OnTheExstnceOfASlwMnfld}
  []{Foias1988_InrtlMnfldsFrNnlnrEvltnryEqtns}
  []{Temam1990_InrtlMnflds}).
%
%
Our aim for designing patch schemes is that the scheme should have both the following properties: have a slow manifold that accurately matches the slow dynamics of the given microscale model; and the scheme should be stable.

\subsection{Extend the staggered grid to multiscale patch scheme}
\label{ssc:extndTheStgrdGrdToMltsclePtchSchme}

To take advantage of the `wave-friendly' features of staggered grids~(\cref{sec:cllctdAndStgrdGrdsFrWvelkeSystms}), here we adapt and extend the collocated patch scheme previously established for dissipative systems \parencite{Roberts2005_HiOrdrAcrcyGapTthScme, Roberts2007_GnrlTthBndryCndnsEqnFreeMdlng, Bunder2017_GdCplingFrTheMltsclePtchSchmeOnSystmsWthMcrscleHtrgnty, Maclean2021_AtlbxOfEqtnFreFnctnsInMtlbOctveFrEfcntSystmLvlSmltn}.
%
\Cref{fig:PtchGrd_n1t0_huvx} shows one example staggered grid of patches, in 2D there are a total of~\(167\,040\) such staggered patch grids possible.
The following list defines the parameters for such staggered patch grids composed of square patches distributed over a square macroscale domain.
\begin{itemize}
  \item \emph{Macro-grid interval}~\(\Delta\), also \emph{inter-patch distance}, is the distance between two adjacent patch centres (size of the {\cIJ violet} grid intervals) in~\(x\)- and \(y\)-directions (we focus on \emph{uniformly spaced patches} over a square domain).
  %
  \item \emph{Number of macro-grid intervals}~\(N\) is the number of grid intervals ({\cIJ violet}) in the periodic domain in each of the~\(x\)- and \(y\)-directions.
  \item \emph{Sub-patch micro-grid interval}~\(\delta\) is the distance between two adjacent micro-grid nodes (size of the {\cij green} grid intervals) in each of the~\(x\)- and \(y\)-directions (we focus on \emph{uniformly spaced micro-grid nodes}).
  \item \emph{Number of sub-patch micro-grid intervals}~\(n\) is the number of microscale grid intervals ({\cij green}) within a square patch in each of the~\(x\)- and \(y\)-directions.
  \item \emph{Patch size}~\(l\) is the side length of the square patch ({\cIJ violet} squares enclosing {\cij green} grid) in~\(x\)- and \(y\)-directions.
  \item \emph{Patch scale ratio}~\(r=l/(2\Delta)\) quantifies the ratio of the simulated to the unsimulated space for each spatial dimension. In practical use, patch scale ratios~\(r = l/(2\Delta)\) are small, typically ranging from~\(0.0001\) to~\(0.1\).
\end{itemize}

\begin{figure}
\caption{\label{fig:allMicroGrids}%
All~\(16\) possible geometrically compatible staggered sub-patch micro-grids (patches) for \pde{}s~\eqref{eqs:stgrdFDEs_gnrcWve} of wave-like systems, with~\(n = 4, 6\) sub-patch micro-grid intervals.
Names indicate the node type on left, right, bottom, and top edges respectively.
Patches with {\cij green} highlighted names have centred nodes.
Centredness depends on both the edge type and whether~\(n/2\) is odd or even.
}
\centering
\begin{subfigure}{0.42\textwidth}
  \centering
  \caption{\label{fig:allMicroGrids_n4}%
  Staggered micro-grids for \(n=4\) (i.e., even~\(n/2\)).
  The patches \(hhhh\), \(uuhh\), and \(hhvv\) are \(h\)-, \(u\)-, and \(v\)-centred patches respectively.
  }
  \includegraphics{figs/allMicroGrids_n4_\nodeShape}
\end{subfigure}%
\quad
\begin{subfigure}{0.5\textwidth}
  \centering
  \caption{\label{fig:allMicroGrids_n6}%
  Staggered micro-grids for \(n=6\) (i.e., odd~\(n/2\)).
  The patches \(uuvv\), \(hhvv\), and \(uuhh\) are \(h\)-, \(u\)-, and \(v\)-centred patches respectively.
  }
  \includegraphics{figs/allMicroGrids_n6_\nodeShape}
\end{subfigure}%
\end{figure}

For a patch scheme over a collocated patch grid, constructing and using a collocated micro-grid (\cref{fig:flDmnGrd_clctd_stgrd}, left) within a patch is straightforward as there is only one way to arrange patch edge nodes.
But when we invoke a staggered micro-grid (\cref{fig:flDmnGrd_clctd_stgrd}, right) within a patch  (the small {\cIJ violet} squares enclosing {\cij green} grid in \cref{fig:PtchGrd_n1t0_huvx}), the heterogeneous nodes lead to many possible arrangements for the patch edge nodes.
The type of edge nodes (\(h\), \(u\), or \(v\)) on left, right, bottom and top edges on a patch edge define the \emph{edge type} of that patch.
We identify a \emph{patch type}, or equivalently a sub-patch \emph{micro-grid type}, by its edge type.
For example, the patch type (or the micro-grid type) \(uuvv\) in the left-bottom of \cref{fig:allMicroGrids_n4,fig:allMicroGrids_n6} has \(u\)-edge nodes on left and right edges, and \(v\)-edge nodes on bottom and top edge nodes.
To simplify the organisation and discussion, micro-grids with \(4\) or \(3\) sub-patch micro-grid intervals are both identified by \(n=4\) as in \cref{fig:allMicroGrids_n4}, similarly micro-grids with \(6\) or \(5\) sub-patch micro-grid intervals are identified by \(n=6\) as in \cref{fig:allMicroGrids_n6}.
But there is no ambiguity as the edge type uniquely identifies them.
When a patch contains a centre node ({\cij green} highlighted names in \cref{fig:allMicroGrids_n4,fig:allMicroGrids_n6}), we call it a \emph{centred patch}, otherwise a \emph{non-centred patch}.
This \emph{centredness} depends on both the patch edge type and whether~\(n/2\) is odd or even.
For example, the \(uuvv\) patch is non-centred (\cref{fig:allMicroGrids_n4}, left-bottom) for \(n = 4\) (i.e., even~\(n/2\)), but the same \(uuvv\) patch is centred (\cref{fig:allMicroGrids_n6}, left-bottom) for \(n = 6\) (i.e., odd~\(n/2\)).
We also identify a centred patch by its centre node.
For example, in \cref{fig:allMicroGrids_n6} for \(n=6\), the \(uuvv\), \(hhvv\), and \(uuhh\) patches are also called \(h\)-, \(u\)-, and \(v\)-patches respectively.
Except where explicitly stated, the details in this article are mainly given for the case of \(n=6\) sub-patch \text{micro-grid intervals.}

\label{pge:dfntn_cmptblePtchGrd}
We define a staggered micro-grid to be \emph{geometrically compatible} if it has all the necessary edge nodes (unfilled markers in \cref{fig:allMicroGrids}) to calculate all the finite differences of the specified \ode{}s~\cref{eqs:stgrdFDEs_gnrcWve} at all the interior nodes.
\Cref{fig:allMicroGrids} shows all~\(16\) possible geometrically compatible staggered sub-patch micro-grids (i.e., patches).
For example, the \(uhvv\) micro-grid (second in the bottom row of \cref{fig:allMicroGrids_n6}) would not be compatible if the right edge all contains \(\vNodeE\,v\)-nodes instead of \(\hNodeE\,h\)-nodes, as the \(\hNodeE\,h\)-nodes on the right edge are necessary to compute~\(\partial h / \partial x\) on the right most interior \(\uNode\,u\)-nodes.

\begin{figure}
\caption{\label{fig:PtchGrd_n1t0_huvx_huvh}%
Example staggered patch grids for \(n=6\) sub-patch micro-grid intervals with three (left) and four patches (right) per macro-cell.
The patch grid parameters are patch size~\(l\), macro-grid interval~\(\Delta\), sub-patch micro-grid interval~\(\delta\).
The patches here are enlarged for visual clarity; in practice the patches are very small relative to their inter-patch distance~\(\Delta\), with small patch scale \text{ratio~\(r = l/(2\Delta)\)} roughly ranging from \(0.0001\) to \(0.1\).
}
\centering
\begin{subfigure}{0.475\textwidth}
  \centering
  \caption{\label{fig:PtchGrd_n1t0_huvx}%
  Patch grid~\#79985, using three patches \(uuvv\), \(hhvv\) and \(uuhh\), gives stable and accurate patch schemes with minimal computation.
  This patch grid is the central focus of this article.
  }
  \input{figs/inkscape/PatchGrid_idlWve_N6n6_n1t0_huvx_\nodeShape.pdf_tex}
\end{subfigure}%
\quad
\begin{subfigure}{0.475\textwidth}
  \centering
  \caption{\label{fig:PtchGrd_n1t0_huvh}%
  Patch grid~\#80001, using four patches \(uuvv\) (twice), \(hhvv\), and \(uuhh\), also gives stable and accurate patch schemes, but with more computation than~\#79985.
  }
  \input{figs/inkscape/PatchGrid_idlWve_N6n6_n1t0_huvh_\nodeShape.pdf_tex}
\end{subfigure}
\end{figure}

Consider a possible 2D staggered patch grid to be designed, comprising the macro-cells ({\cpq orange} squares in \cref{fig:PtchGrd_n1t0_huvx}) each containing \(2 \times 2 = 4\) possible patches ({\cIJ violet} squares enclosing {\cij green} micro-grid in \cref{fig:PtchGrd_n1t0_huvx}).
These macro-cells in a patch grid are on the macroscale similar to the micro-cells in a micro-grid ({\cpq orange} square in \cref{fig:flDmnGrd_clctd_stgrd}) on the microscale.
Similar to the \(n^2/4\)~micro-cells in a \(n \times n\) full-domain micro-grid\footnote{We use the same symbol~\(n\) and~\(\delta\) for both the full-domain micro-grid and the sub-patch micro-grid, and disambiguate by words and/or context.} (\cref{sec:cllctdAndStgrdGrdsFrWvelkeSystms}), there are \(N^2/4\)~macro-cells in a \(N \times N\) staggered patch grid.
Each of the four patches within a macro-cell could be either empty or contain one of the \(16\)~possible micro-grids of \cref{fig:allMicroGrids} depending upon whether~\(n/2\) is odd or even.
Thus, excluding the all-empty case, in a \(2 \times 2\) macro-cell configuration, the total number of possible 2D staggered patch grids is \((16+1)^4 - 1 = 83\,520\) for each~\(n\).
Including all the \(32\)~possible sub-patch micro-grids in \cref{fig:allMicroGrids} for any number of sub-patch micro-grid intervals~\(n\), there are a total of~\(167\,040\) possible 2D staggered patch grids.
We analysed all of these~\(167\,040\) staggered patch grids.
For each of the two cases of~\(n/2\) being odd or even, this article separately analyses the respective \(83\,520\) staggered patch grids, covering the total of~\(167\,040\) patch grids.
We uniquely identify each of these \(83\,520\) staggered patch grids for a particular~\(n\) using a base~\(17\) number whose decimal equivalent, \emph{patch grid Id}, ranges from~\#1 to~\#83520.
For a particular~\(n\), each of the \(17\)~digits encodes either an empty patch or one of the \(16\)~possible sub-patch micro-grids in \cref{fig:allMicroGrids}.
The magnitude of these numbers does not have any significance, but the closer these numbers the more similar the patch grid geometries are.
For example, \cref{fig:PtchGrd_n1t0_huvx_huvh} shows two staggered patch grids with Ids~\#79985 and \#80001; these two patch grids are nearly the same except for an additional \(h\)-centred patch in each macro-cell for the \#80001 patch grid.

Based on the centredness of the sub-patch micro-grids ({\cij green} highlighted in \cref{fig:allMicroGrids}), there are two types of qualitatively different staggered patch grids:
\begin{itemize}
  \item \emph{centred staggered patch grids}\label{pge:dfntn_cntrdPtchGrid} whose sub-patch micro-grids all have a centre node and include a centred-patch corresponding to each of the state variables;
  %
  \item staggered patch grids that are not centred staggered patch grids are \emph{non-centred staggered patch grids}.
\end{itemize}
For example, for \(n=6\) the two patch grids in \cref{fig:PtchGrd_n1t0_huvx_huvh}, \#79985 and \#80001, are both centred patch grids.
%
%
But in \cref{fig:PtchGrd_n1t0_huvx_huvh} replacing any of the centred patch or the empty patch with a non-centre patch gives a non-centred patch grid.
Even when all sub-patch micro-grids have a centre node, if a patch grid does not contain a centred-patch corresponding to each of the state variables, it is not a centred patch grid.
For example, replacing all the centred patches in \cref{fig:PtchGrd_n1t0_huvx} with any one particular centred patch, say all-\(h\)-centred patches, gives a non-centred patch grid.
%
For patches without a centre node, we compute the macroscale patch value by averaging over the closest respective~\(h, u, v\) values to the patch centre.
%
It is these macroscale patch values that are used in coupling the patches (\cref{ssc:ptchCplngCnctsTheScls}).

In \cref{fig:PtchGrd_n1t0_huvx} with three non-empty patches, any of the \(24\)~permutations of the \(h\)-, \(u\)-, \(v\)-patches or an empty patch is also a centred patch grid as per the definition.
Hence, for a given~\(n\) there are \(24\)~centred staggered patch grids containing three patches per cell.
In \cref{fig:PtchGrd_n1t0_huvh} with four non-empty patches, any of the \(12\)~unique permutations of the \(h\)-, \(u\)-, \(v\)-patches with repeated \(h\)-patch is also a centred patch grid.
%
A similar patch grid with a repeated \(u\)-patch or a repeated \(v\)-patch has \(12\)~unique permutations of centred patch grids.
Hence, for a given~\(n\) there are \(36\)~centred staggered patch grids containing four patches per cell.
%
So, for a given~\(n\), among the \(83\,520\) staggered patch grids, \(24+36 = 60\)~of them are centred staggered patch grids containing either three or four patches per cell; and the remaining~\(83\,460\) are non-centred patch grids.
%
Thus, for both the cases of~\(n/2\) being odd or even, in total there are \(120\)~centred staggered patch grids and \(166\,920\) non-centred patch grids.

Whereas all the possible~\(167\,040\) staggered patch grids are geometrically compatible 2D discretisations for simulating multiscale wave physics, most of them constitute unstable patch schemes (\cref{ssc:onlyPtchGrdsWthAlSymtrcPtchsAreStble}).
The geometry-stability study of \cref{ssc:onlyPtchGrdsWthAlSymtrcPtchsAreStble} shows that \emph{almost all non-centred staggered patch grids lead to unstable patch schemes but all centred patch grids constitute stable patch schemes}
The accuracy study of \cref{ssc:onlyCntrdPtchGrdsAreAcrte} shows that all the centred staggered patch grids constitute accurate patch schemes and that none of the stable non-centred patch grids constitutes accurate patch schemes.
That is, only the \(120\)~centred staggered patch grids among all the \(167\,040\) compatible staggered patch grids constitute stable and accurate multiscale patch schemes.
All \(120\)~centred staggered patch grids that constitute stable patch schemes are with the same accuracy (\cref{ssc:onlyCntrdPtchGrdsAreAcrte}), yet due to fewer nodes, the \(48\)~centred patch grids with three patches per macro-cell cost less computation compared to the remaining \(72\)~centred patch grids with four patches per macro-cell.

The staggered patch grid~\#79985 in \cref{fig:PtchGrd_n1t0_huvx} is the main example for the specific discussions in this article.

\begin{figure}
  \caption{\label{fig:indexConvention}%
  Index convention, using the patch grid~\#79985 in \cref{fig:PtchGrd_n1t0_huvx} as a specific example, but the same convention is used for all the \(167\,040\) compatible 2D staggered patch grids.
  Here, the staggered patch grid has~\(N=6\) macro-grid intervals ({\cIJ violet}), containing staggered patches of staggered micro-grids with~\(n=6\) micro-grid intervals ({\cij green}).
  }
  \centering
  \begin{subfigure}[b]{0.47\textwidth}
    \centering
    \caption{\label{fig:indexing_sbPtchMicroGrd}%
      Sub-patch micro-grid node index~\({\cij i,j} \in \{1,2,\ldots,n-1\}\) for the interior nodes (\(\hNode\), \(\uNode\), \(\vNode\)) and \({\cij i,j} \in \{0, n\}\) for patch edge nodes (\(\hNodeE\), \(\uNodeE\), \(\vNodeE\)).
      }
    \footnotesize
    \input{figs/inkscape/indexing_subPatchMicroGrid_\nodeShape.pdf_tex}
  \end{subfigure}%
  \quad
  \begin{subfigure}[b]{0.47\textwidth}
    \centering
    \caption{\label{fig:indexing_patchGrid}%
      Global (macroscale) patch index~\({\cIJ I,J} \in \{0,1,\ldots,N-1\}\);
      local sub-cell patch index  \({\cpq p,q} \in \{0,1\}\) within each macro-cell ({\cpq orange} squares).}
    \footnotesize
    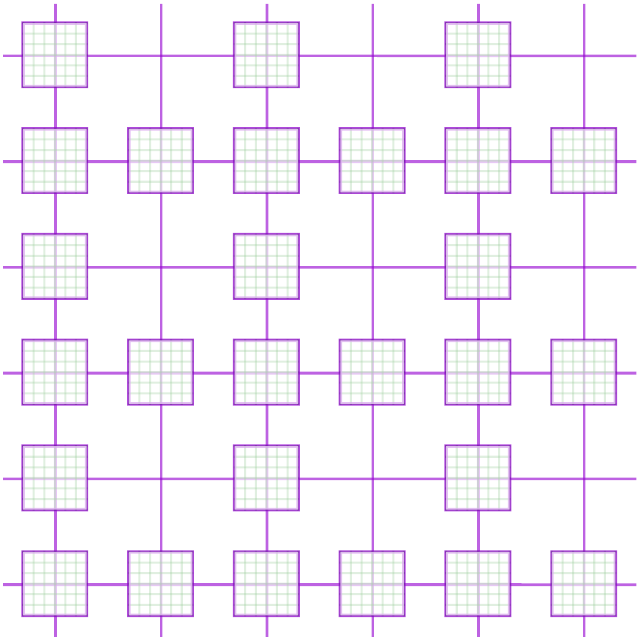
  \end{subfigure}
\end{figure}

\needspace{3\baselineskip}
\Cref{fig:indexConvention} illustrates the following three kinds of indices we use to identify patches and sub-patch micro-grid nodes in all cases, using the staggered patch grid~\#79985 in \cref{fig:PtchGrd_n1t0_huvx} as a specific example:
\begin{itemize}
  \item the pair \({\cIJ I,J} \in \{0, 1, \ldots, N-1\}\) is the global (macroscale) patch index;

  \item the pair \({\cpq p,q}\), defined by \({\cpq p} := {\cIJ I}\bmod 2\) and \({\cpq q} := {\cIJ J}\bmod 2\)\,, is the local (macroscale) sub-cell patch index, that is, \({\cpq p,q} \in \{0,1\}\) within each macro-cell ({\cpq orange} squares in \cref{fig:indexing_patchGrid});

  \item the pair \({\cij i,j}\) is the sub-patch micro-grid node index, with~\({\cij i,j} \in \{1,\ldots,n-1\}\) for interior nodes (filled markers~\(\hNode\), \(\uNode\), \(\vNode\) in \cref{fig:indexing_sbPtchMicroGrd}) and \({\cij i,j} \in \{0,n\}\) for patch edge nodes (unfilled markers~\(\hNodeE\), \(\uNodeE\), \(\vNodeE\) in \cref{fig:indexing_sbPtchMicroGrd}).
\end{itemize}

Using the microscale model~\cref{eqs:stgrdFDEs_gnrcWve} within the patches in a staggered patch grid (\cref{fig:PtchGrd_n1t0_huvx}) with~\(N\) macro-grid intervals (i.e., total of \(N/2 \times N/2\) macro-cells), where each patch consists of \(n \times n\) sub-patch micro-grid intervals (e.g., \(n=6\) in \cref{fig:PtchGrd_n1t0_huvx}), gives the \emph{staggered patch scheme} for the microscale model~\cref{eqs:stgrdFDEs_gnrcWve} of a generic wave-like system as the following in the specific case of \cref{fig:PtchGrd_n1t0_huvx} (with \({\cij i,j} \in \{1,\ldots,n-1\}\)):
\begin{subequations}  \label{eqs:FDEs_pN_gnrcWve}
\allowdisplaybreaks
\begin{alignat}{2}
   \label{eqn:FDEs_pN_gnrcWve_h}
   \hNode\;
    &\frac{d}{d t} h^{{\cIJ I,J}}_{{\cij i,j}}(t) &&=
    -\frac{u^{{\cIJ I,J}}_{\cij i+1,j} - u^{\cIJ I,J}_{\cij i-1,j}}{2 \delta}
      - \frac{v^{\cIJ I,J}_{\cij i,j+1} - v^{\cIJ I,J}_{\cij i,j-1}}{2 \delta}
      + \cdots
      \,,\\ \notag &&&\qquad
       {\cij i}-{\cpq p}\text{ odd,}\quad
       {\cij j}-{\cpq q}\text{ odd};
      \\
  \label{eqn:FDEs_pN_gnrcWve_u}
  \uNode\;
    &\frac{d}{d t} u^{{\cIJ I,J}}_{{\cij i,j}}(t) &&=
    - \frac{
      h^{\cIJ I,J}_{i + 1,j}
      - h^{\cIJ I,J}_{i - 1,j}
      }{
      2 \delta
      }
    + \cdots \,,\quad
       {\cij i}-{\cpq p}\text{ even,}\quad
       {\cij j}-{\cpq q}\text{ odd};
      \\
  \label{eqn:FDEs_pN_gnrcWve_v}
  \vNode\;
    &\frac{d}{d t} v^{{\cIJ I,J}}_{{\cij i,j}}(t) &&=
    - \frac{
        h^{\cIJ I,J}_{i,j + 1}
        - h^{\cIJ I,J}_{i,j - 1}
        }{
        2 \delta
        }
    + \cdots
       \,,\quad
       {\cij i}-{\cpq p}\text{ odd,}\quad
       {\cij j}-{\cpq q}\text{ even};
\end{alignat}
\emph{and a patch coupling} computes the edge values \(\hNodeE\,\hVarP\), \(\uNodeE\,\uVarP\), and \(\vNodeE\,\vVarP\) for
\begin{align}
\parbox{25em}{\raggedright\let\quad\ 
    \({\cij i} \in \{0, n\}\) for left and right edges,  \quad   
    \({\cij j} \in \{0, n\}\) for bottom and top edges, \quad
  \({\cIJ I}\) even for~\(h,v\),
  \quad  \({\cIJ I}\) odd for~\(u\),
  \quad \({\cIJ J}\) even for~\(h,u\),
  \quad \({\cIJ J}\) odd for~\(v\).
}
\end{align}
\end{subequations}
Due to the chosen periodic boundary conditions for the full-domain microscale model~\cref{eqs:stgrdFDEs_gnrcWve}, the three fields~\(h\), \(u\), \(v\) are macroscale \(N\)-periodic in both~\(\cIJ I\) and~\(\cIJ J\), where \(N = L/\Delta\).
A specific patch coupling method computes patch edge nodes (\(\hNodeE\,\hVarP\), \(\uNodeE\,\uVarP\), \(\vNodeE\,\vVarP\) in \cref{fig:PtchGrd_n1t0_huvx}) from the aggregate values of the neighbouring patches
(\(\hNode\,\hVarP\), \(\uNode\,\uVarP\), \(\vNode\,\vVarP\) with~\({\cij i} = {\cij j} = n/2\) in the case of \cref{fig:PtchGrd_n1t0_huvx}).
The patch coupling provides the mechanism whereby patches influence each other.

Arranging the patch interior values of the scheme~\cref{eqs:FDEs_pN_gnrcWve}, which are over the staggered patch grid such as in \cref{fig:indexConvention}, into a vector gives the state vector~\(\vec{x}^{I}\) of the staggered patch scheme, which is a dynamic state variable evolving in time.
Here the superscript~\({(\argdot)^I}\) is not an index or exponent, instead, a qualifier denoting the patch interior nodes. 
In the specific case of \cref{fig:indexConvention}, the staggered patch scheme state vector
\begin{equation*} 
  \vec{x}^{I} = \parbox[t]{0.9\linewidth}{\mathcode`\,="213B$
  \big(
  h_{{\cij{1,1}}}^{\cIJ{0,0}},
  h_{{\cij{1,3}}}^{\cIJ{0,0}}
  ,\ldots,
  u_{{\cij{2,1}}}^{\cIJ{0,0}},
  u_{{\cij{2,3}}}^{\cIJ{0,0}}
  ,\ldots,
  v_{{\cij{1,2}}}^{\cIJ{0,0}},
  v_{{\cij{1,4}}}^{\cIJ{0,0}}
  ,\ldots,
  h_{{\cij{1,2}}}^{\cIJ{0,1}},
  h_{{\cij{1,4}}}^{\cIJ{0,1}}
  ,\ldots,
  u_{{\cij{2,2}}}^{\cIJ{0,1}},
  u_{{\cij{2,4}}}^{\cIJ{0,1}}
  ,\ldots,
  v_{{\cij{1,1}}}^{\cIJ{0,1}},
  v_{{\cij{1,3}}}^{\cIJ{0,1}}
  ,\ldots,
  h_{{\cij{1,1}}}^{\cIJ{0,2}},
  h_{{\cij{1,3}}}^{\cIJ{0,2}}
  ,\ldots,
  u_{{\cij{2,1}}}^{\cIJ{0,2}},
  u_{{\cij{2,3}}}^{\cIJ{0,2}}
  ,\ldots,
  v_{{\cij{1,2}}}^{\cIJ{0,2}},
  v_{{\cij{1,4}}}^{\cIJ{0,2}}
  ,\ldots, h_{{\cij{2,1}}}^{\cIJ{1,0}},
  h_{{\cij{2,3}}}^{\cIJ{1,0}}
  ,\ldots,
  u_{{\cij{1,1}}}^{\cIJ{1,0}},
  u_{{\cij{1,3}}}^{\cIJ{1,0}}
  ,\ldots,
  v_{{\cij{2,2}}}^{\cIJ{1,0}},
  v_{{\cij{2,4}}}^{\cIJ{1,0}}
  ,\ldots,
  h_{{\cij{2,1}}}^{\cIJ{1,2}},
  h_{{\cij{2,3}}}^{\cIJ{1,2}}
  ,\ldots,
  u_{{\cij{1,1}}}^{\cIJ{1,2}},
  u_{{\cij{1,3}}}^{\cIJ{1,2}}
  ,\ldots,
  v_{{\cij{2,2}}}^{\cIJ{1,2}},
  v_{{\cij{2,4}}}^{\cIJ{1,2}}
  ,\ldots,
  h_{{\cij{1,1}}}^{\cIJ{0,2}},
  h_{{\cij{1,3}}}^{\cIJ{0,2}}
  ,\ldots,
  u_{{\cij{2,1}}}^{\cIJ{0,2}},
  u_{{\cij{2,3}}}^{\cIJ{0,2}}
  ,\ldots,
  v_{{\cij{1,2}}}^{\cIJ{0,2}},
  v_{{\cij{1,4}}}^{\cIJ{0,2}}
  ,\ldots \big) .
$}
\end{equation*}
Here, the total number of patch interior nodes, the size of the state vector~\(\vec{x}^I\),
\begin{equation} \label{eqn:exprsn_nI_gLnrWve}
  n_p^I := (N^2/4) (9n^2/4 - 4n + 2)\,,
\end{equation}
where \(N\)~is the number of macro-grid intervals and \(n\)~is the number of sub-patch micro-grid intervals.
For example, for \(N = 6\), \(10\), \(14\), \(18\), \(22\), \(26\) macro-grid intervals with~\(n=6\) sub-patch micro-grid intervals, \(n_p^I = 531\), \(1475\), \(2891\), \(4779\), \(7139\), \(9971\) respectively.
This~\(n_p^I\) is only for the \(48\) centred patch grids with three patches per macro-cell, the \(72\) centred patch grids with four patches per macro-cell have more interior nodes.

Arranging the patch edge values of the staggered patch scheme~\cref{eqs:FDEs_pN_gnrcWve}, with~\({\cij i,j} \in \{0,n\}\) for \cref{fig:indexConvention}, into a vector gives the edge vector~\(\vec{x}^{E}\), of size~\(n^E_p\), containing all the edge values of all the patches.
The edge vector~\(\vec{x}^{E}\) is a function of the state vector~\(\vec x^I\)---a function~\(\vec{x}^E(\vec{x}^{I})\) determined by the patch coupling of a particular \text{patch scheme.}

Using a multiscale staggered patch grid with one layer of edge nodes, such as in~\cref{fig:PtchGrd_n1t0_huvx_huvh}, allows for calculating the first spatial derivatives on all patch interior nodes.
For ideal waves modelled by \pde{}s~\cref{eqs:PDEs_gnrcWve}, after dropping the additional terms denoted by ``\(\cdots\)'', a staggered patch grid with one layer of edge nodes is sufficient as in~\cref{fig:PtchGrd_n1t0_huvx_huvh,fig:indexing_sbPtchMicroGrd}.
The staggered sub-patch micro-grids in~\cref{fig:PtchGrd_n1t0_huvx_huvh} are also sufficient for many other wave \pde{}s containing only first-order spatial derivatives.
Multiscale staggered grid discretisation of more complex terms in the \pde{}s~\cref{eqs:PDEs_gnrcWve} often requires additional layers of edge nodes (\cref{ssc:mltiLyrEdgeNdsFrHghrOrdrSptlDrvtvs}).
For the specific multiscale staggered patch grid in \cref{fig:PtchGrd_n1t0_huvx}, the total number of patch edge nodes, the size of the edge vector~\(\vec{x}^E\), is
%
  \(n^E_p := (N^2/4) (5n - 4)\)\,,
%
where \(N\)~is the number of macro-grid intervals and \(n\)~is the number of sub-patch micro-grid intervals.
For example, for the staggered patch grid in \cref{fig:PtchGrd_n1t0_huvx} with \(N = 6\), \(10\), \(14\), \(18\), \(22\), \(26\), and  \(n=6\), \(n^E_p = 360\), \(1000\), \(1960\), \(3240\), \(4840\), \(6760\) respectively.
The \(72\)~centred patch grids with four patches per macro-cell have more edge nodes.

In terms of the state vector~\(\vec{x}^{I}\) and the edge vector function~\(\vec{x}^{E}(\vec{x}^{I})\),  a staggered patch scheme~\cref{eqs:FDEs_pN_gnrcWve} is represented as a dynamical system by \ode{}s of the form
\begin{equation} \label{eqn:dynSys_pN_gnrcWve}
  \dby{\vec{x}^I}{t} = \vec{F}\big( \vec{x}^{I};\, \vec{x}^{E}(\vec{x}^{I}) \big).
\end{equation}
The \(\vec{F}\big( \vec{x}^{I}; \vec{x}^{E}(\vec{x}^{I}) \big)\) in the \ode{}s~\cref{eqn:dynSys_pN_gnrcWve} corresponds to the~\(\vec{f}( \vec{x} )\) in the  microscale model~\cref{eqs:dynSys_mdN_gnrcWve}.
The functions~\(\vec{F}\) and~\(\vec{f}\) encode the same microscale model for the generic wave-like system~\cref{eqs:PDEs_gnrcWve}, except for the following difference:
the only argument of~\(\vec{f}\) is the nodal values~\(\vec{x}\) in the full-domain micro-grid;
whereas the two arguments of~\(\vec{F}\) are the patch interior values~\(\vec{x}^I\) and the patch edge values~\(\vec{x}^E(\vec x^I)\) from the patch coupling.

A staggered patch scheme simulation is performed by numerical time-integration of the \ode{}s~\cref{eqn:dynSys_pN_gnrcWve} for the particular staggered microscale \ode{}s~\cref{eqs:FDEs_pN_gnrcWve} on the interior nodes of the patches (e.g., filled markers in \cref{fig:PtchGrd_n1t0_huvx,fig:mltiLyrEdgeNdsForPatchGrid}), with some particular patch coupling represented by~\(\vec{x}^{E}(\vec{x}^{I})\).

\subsection{Multi-layer edge nodes for higher-order spatial derivatives}
\label{ssc:mltiLyrEdgeNdsFrHghrOrdrSptlDrvtvs}

The staggered patch grids of \cref{fig:PtchGrd_n1t0_huvx_huvh} are suitable for \pde{}s with only first order spatial derivatives, such as \pde{}s~\cref{eqs:PDEs_gnrcWve} with ``\(+\cdots\)'' omitted.
However, sub-patch microscale discretisations of higher-order spatial derivatives (such as second spatial derivative for viscous diffusion) need to use a wider stencil of nodes on the micro-grid.
A wider stencil is not an issue for most of the interior nodes, but the interior nodes closest to the edges of a patch need values of nodes that lie outside the sub-patch micro-grids in \cref{fig:PtchGrd_n1t0_huvx_huvh}.

As a prototypical example of wave-like \pde{}s~\cref{eqs:PDEs_gnrcWve} with high-order derivatives, we consider the non-dimensional \text{\emph{weakly damped linear wave \pde{}s}}
\begin{subequations} \label{eqs:PDEs_gLnrWve}
\begin{align}
  \label{eqn:PDE_gLnrWve_h}
  \doh{h}{t} &= - \doh{u}{x} - \doh v y\,,\\
  \label{eqn:PDE_gLnrWve_u}
  \doh{u}{t} &= - \doh{h}{x} - c_D u + c_V \doh[2]{u}{x} + c_V \doh[2]{u}{y} \,, \\
  \label{eqn:PDE_gLnrWve_v}
  \doh{v}{t} &= - \doh{h}{y} - c_D v \, + c_V \doh[2]{v}{x} + c_V \doh[2]{v}{y}\,,
\end{align}
\end{subequations}
with linear drag and viscous diffusion respectively characterised by the coefficients~\(c_D, c_V\).
The macroscale boundary conditions remain that the three fields~\(h\), \(u\), and~\(v\) are \(L\)-periodic in both~\(x\) and~\(y\).
The case \(c_D = c_V = 0\) corresponds to the non-dissipative \emph{ideal wave \pde{}s} (\cites[pp.136--137]{Dean1991_WtrWveMchncsFrEngnrsAndScntsts}[pp.257--258]{Mehaute1976_AnIntro2HydrdynmcsWtrWvs}). 
General solutions of the \pde{}s~\eqref{eqs:PDEs_gLnrWve} are linear combinations of progressive sinusoidal waves that decay at rates depending upon~\(c_D\) and~\(c_V\).

We explored the following two families of methods of approximating the second derivatives in~\eqref{eqs:PDEs_gLnrWve} on the micro-grids inside each patch.
\begin{itemize}
  \item This first methods was unsuccessful.
  Whenever the microscale finite difference model needs nodes outside the patch edges, we used the nearest edge and interior nodes to extrapolate the required outside node values (we tried constant values, linear and quadratic extrapolation).
  In effect, such extrapolations change the centred finite differences, such as in~\cref{eqs:FDEs_pN_gnrcWve}, into non-centred finite differences for the interior nodes close to the edge nodes.
  But it eventuated that such extrapolation makes most variants of the staggered patch schemes unstable.

\begin{figure}
\caption{\label{fig:mltiLyrEdgeNdsForPatchGrid}%
  Centred staggered patch grid~\#79985 of \cref{fig:PtchGrd_n1t0_huvx} with the same grid parameters but with additional multi-layer edge nodes, contain the same number of interior nodes \(\hNode\,\hVarP\), \(\uNode\,\uVarP\), \(\vNode\,\vVarP\) for \({\cij i,j} \in \{1,\ldots,n-1\}\).
  But the number of edge values \(\hNodeE\,\hVarP\), \(\uNodeE\,\uVarP\), and \(\vNodeE\,\vVarP\) depend on the number of layers of edge nodes.
  Layer numbers are shown for the right edge.
}
\begin{subfigure}[b]{0.49\textwidth}
  \centering
  \caption{\label{fig:PtchGrd_n2t0}%
    Two layers of edge nodes in normal direction to the edges, no edge nodes in tangential direction to the edges (e.g., no \(\vNode\,v\) node on corners of \(\uNode\,u\)-centred patch).
  }
  \input{figs/inkscape/PatchGrid_gLnrWve_N6n6_n2t0_\nodeShape.pdf_tex}
\end{subfigure}
\;
\begin{subfigure}[b]{0.49\textwidth}
  \centering
  \caption{\label{fig:PtchGrd_n3t2}%
    Three layers of edge nodes in normal direction to the edges, two layers of edge nodes in tangential direction to the edges.
  }
  \input{figs/inkscape/PatchGrid_vSwtrC_N6n6_n3t2_\nodeShape.pdf_tex}
\end{subfigure}
\end{figure}

  \item A successful alternative method that we developed is to append additional layers of edge nodes to the sub-patch micro-grid as in \cref{fig:PtchGrd_n2t0,fig:PtchGrd_n3t2}.
  We calculate the values of these extra edge layer nodes using the same spectral interpolation as in calculating the edge values of \cref{fig:PtchGrd_n1t0_huvx}.
  Unlike the first method using extrapolation, this method uses the same spatial finite difference formula at all interior nodes within each patch, including those close to the edge nodes.
  \cref{sec:mltscleStgrdPtchGrdFrWklyDmpdLnrWvs} shows this method of additional layers of edge nodes gives good families of stable and accurate patch schemes.
\end{itemize}

\Cref{fig:mltiLyrEdgeNdsForPatchGrid} shows two example staggered patch grids with additional layers of edge nodes.
%
Irrespective of the number of layers of the edge nodes, the staggered patch grids contain the same number of patch interior nodes for a given number of sub-patch micro-grid intervals (e.g., \(n=6\) in \cref{fig:PtchGrd_n1t0_huvx_huvh,fig:mltiLyrEdgeNdsForPatchGrid}).
For finite difference approximation of all the terms in the \pde{}s~\cref{eqs:PDEs_gLnrWve}, a staggered patch grid in \cref{fig:PtchGrd_n2t0} is sufficient (\cref{ssc:drveEgnvlsOfAMltscleStgrdPtchSchmeFrLnrMicrscleMdls}).
For \cref{fig:PtchGrd_n2t0} with two layers of edge nodes in the normal direction to the edges and no edge nodes in the tangential direction to the edges, the left and right edge values are \(\hNodeE\,\hVarP\), \(\uNodeE\,\uVarP\), \(\vNodeE\,\vVarP\), for \({\cij i} \in \{-1,0,\; n,n+1\}\) and \({\cij j} \in \{1, 2, \ldots, n-1\}\).
Similarly the bottom and top edge value indices are \({\cij i} \in \{1, 2, \ldots, n-1\}\) and \({\cij j} \in \{-1,0,\; n,n+1\}\).

Nonlinear \pde{}s with higher-order and/or mixed derivatives may require patches with more edge node layers as depicted in \cref{fig:PtchGrd_n3t2}.
The total number of edge nodes depends on the number of layers of the edge nodes.
Our subsequent articles explore such staggered patch grids with more edge node layers for nonlinear wave \pde{}s, specifically for the viscous shallow water \pde{}s of \textcite{Roberts2006_AnAccrteAndCmprhnsveMdlOfThnFldFlwsWthInrtaOnCrvdSbstrts} and the turbulent shallow water \pde{}s of \textcite{Cao2016_MdlngSspnddSdmntInEnvrnmntlTrblntFlds}.


\section{Multiscale staggered patch grid for weakly damped linear waves}
\label{sec:mltscleStgrdPtchGrdFrWklyDmpdLnrWvs}

This section shows that all the \(60\)~centred multiscale staggered patch grids (defined in \cref{ssc:extndTheStgrdGrdToMltsclePtchSchme}) give stable and accurate patch schemes for weakly damped linear wave \pde{}s~\cref{eqs:PDEs_gLnrWve} including the ideal waves with~\(c_D = c_V = 0\).
This section focuses on the \emph{spectral} patch scheme over 2D staggered patch grids for a representative set of physical parameters~\(c_D\), \(c_V\), and discretisation parameters~\(n\), \(\delta\), and \(\Delta\).
A subsequent article will establish the stability, accuracy and sensitivity of several patch coupling schemes for a wide range of parameters.

\Cref{ssc:onlyPtchGrdsWthAlSymtrcPtchsAreStble} studies the dependence of the patch scheme stability on the patch grid geometry, the \emph{geometry-stability study}, for all the \(167\,040\) compatible 2D staggered patch grids (defined in \cref{ssc:extndTheStgrdGrdToMltsclePtchSchme}).
The geometry-stability study of \cref{ssc:onlyPtchGrdsWthAlSymtrcPtchsAreStble} shows that among all the compatible~\(167\,040\) compatible patch grids, only \(1248\) patch grids (\(0.75\%\)) are stable---they all have symmetric patches.
%
Based on some examples of unstable staggered patch grids, \cref{ssc:onlyPtchGrdsWthAlSymtrcPtchsAreStble} also shows that the more non-centred patches in a staggered patch grid, the more unstable modes in the patch scheme.

\Cref{ssc:onlyCntrdPtchGrdsAreAcrte} shows that none of the non-centred patch grids that are stable (\cref{ssc:onlyPtchGrdsWthAlSymtrcPtchsAreStble}) is accurate.
Thus, among all the possible~\(167\,040\) compatible 2D staggered patch grids, only \(120\)~centred patch grids (\(0.07\%\)) are both stable and accurate.

Recall the patches are distributed on a macroscale spatial grid that is translationally invariant to macroscale shifts by multiples of~\(2\Delta\).
Hence the macroscale, small wavenumber, modes of the patch scheme have a purely sinusoidal spatial structure---\emph{exactly the same} sinusoidal spatial structure as that of the small wavenumber modes of the underlying microscale discrete model~\eqref{eqs:stgrdFDEs_gnrcWve}.
Hence, in this scenario, the only error of the \emph{macroscale} waves in the patch scheme is in their time dependence, which is completely characterised by their eigenvalues.
Consequently, we establish the stability and accuracy of the patch scheme by comparing its \emph{macroscale} eigenvalues with the corresponding eigenvalues of the full-domain microscale discretisation~\eqref{eqs:stgrdFDEs_gnrcWve} for weakly damped linear waves \textsc{pde}s~\eqref{eqs:PDEs_gLnrWve}.
\Cref{ssc:drveEgnvlsOfAMltscleStgrdPtchSchmeFrLnrMicrscleMdls} explains the method of deriving eigenvalues of a multiscale patch scheme~\cref{eqs:FDEs_pN_gnrcWve} on a staggered patch grid.

\begin{figure}
\centering
  \begin{minipage}[t][][c]{0.58\linewidth}
    \caption{\label{fig:cmplxEigs_idlWve_Spectral_N10_n6_r0p01_huvx_oce}%
    Eigenvalues for non-dissipative ideal wave \pde{}s~\cref{eqs:PDEs_gLnrWve} with~\(c_D = c_V = 0\), over a \(2\pi \times 2\pi\) periodic domain.
    Eigenvalues~\(\lambda_p\) are for the spectral patch scheme~\cref{eqs:FDEs_pN_gnrcWve} over the centred staggered patch grid~\#79985 (\cref{fig:PtchGrd_n1t0_huvx}) with \(N=10\), \(n=6\), \(r=0.01\).
    Eigenvalues~\(\lambda_\mu\) are for the staggered full-domain microscale model~\eqref{eqs:stgrdFDEs_gnrcWve}.
    The micro-grid interval~\(\delta = 2\pi/3000\) for both.
    The pure imaginary eigenvalues~\(\lambda_p\) show that our patch scheme preserves the wave nature of the underlying system.
    The slow macroscale eigenvalues of the patch scheme (\(\left|\Im(\lambda_p)\right| < 10\)) agree to graphical accuracy with the corresponding slow eigenvalues~\(\lambda_\mu\) of the full-domain microscale model.
    }
    \end{minipage}
    \quad 
      \raisebox{-\height}{%
      \includegraphics{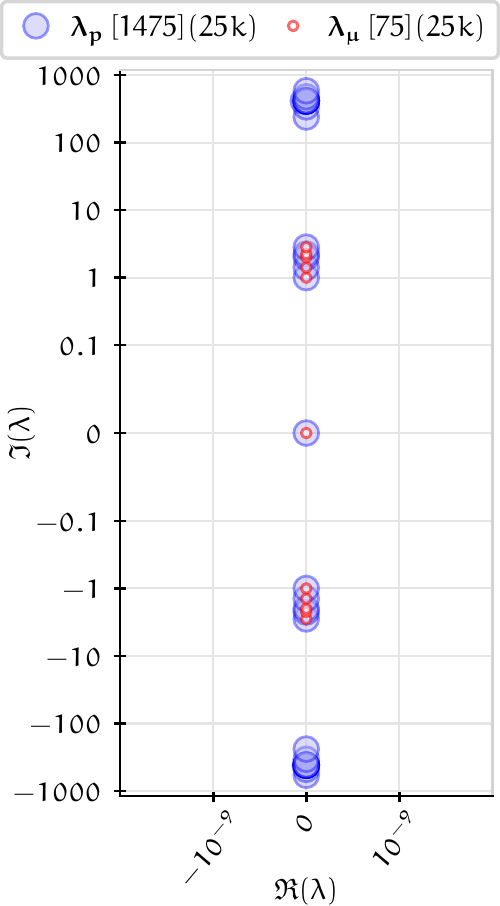}
      }
\end{figure}

For a first example, \cref{fig:cmplxEigs_idlWve_Spectral_N10_n6_r0p01_huvx_oce} shows the eigenvalues of the spectral patch scheme (\cref{ssc:drveEgnvlsOfAMltscleStgrdPtchSchmeFrLnrMicrscleMdls}) over the  multiscale staggered patch grid~\#79985 of \cref{fig:PtchGrd_n1t0_huvx} applied to the non-dissipative ideal waves~\eqref{eqs:FDEs_pN_gnrcWve} (when ``\(+\cdots\)'' are omitted).
As in \cref{fig:cmplxEigs_idlWve_Spectral_N10_n6_r0p01_huvx_oce}, all the eigenvalue plots in this article have a legend that within brackets,~\([\argdot]\), give the total number of eigenvalues for each set.
Also, the legends give the total number of wavenumbers that each set of eigenvalues corresponds to within parentheses,~\((\argdot k)\); the~\(k\) here denotes the wavenumber~\(k\), not thousands.
That the eigenvalues~\(\lambda_p\) are all pure imaginary shows that the patch grid preserves the wave nature of the underlying microscale model, without introducing any significant artificial dissipation.
The slow macroscale eigenvalues of the patch scheme (\(\left|\Im(\lambda_p)\right| < 10\)) agree to graphical accuracy with the corresponding slow eigenvalues~\(\lambda_\mu\) of the full-domain microscale model~\eqref{eqs:stgrdFDEs_gnrcWve}.
\Cref{tbl:macroEigsAgreeExactly} of \cref{ssc:onlyCntrdPtchGrdsAreAcrte} shows this agreement is exact (within~\(10^{-12}\)).
This agreement implies that the slow macroscale dynamics on this staggered patch grid in \cref{fig:PtchGrd_n1t0_huvx} accurately model the corresponding slow macroscale dynamics of the underlying microscale model.

As in \cref{fig:cmplxEigs_idlWve_Spectral_N10_n6_r0p01_huvx_oce}, all the eigenvalue plots in this article are on a quasi-log nonlinear scale using the \(\arcsinh\) function.
%
Specifically, we plot eigenvalues~\(\lambda\) as points on the complex plane with real part~\(\arcsinh(10^{6} \cdot \Re \lambda)\) and complex part~\(\arcsinh( 100 \cdot \Im \lambda)\).
This nonlinear \(\arcsinh\) scale empowers us to see very small as well as very large values all together in the same plot, but unlike a log scale, for both positive and negative values.

\needspace{4\baselineskip}
\subsection{Derive eigenvalues of a multiscale staggered patch scheme for linear microscale models}
\label{ssc:drveEgnvlsOfAMltscleStgrdPtchSchmeFrLnrMicrscleMdls}

%
The method of deriving patch scheme eigenvalues in this subsection is applicable for linear wave-like \pde{}s over any of the \(167\,040\) compatible 2D staggered patch grids, irrespective of the number of edge node layers in the patch scheme (\cref{sec:mltscleStgrdPtchGrdFrWklyDmpdLnrWvs}).
%
But the details of the discussions are for the specific staggered patch grid~\#79985 of \cref{fig:PtchGrd_n2t0} with two layers of edge nodes in the normal direction to the edges, no edge nodes in the tangential direction to the edges.

First, to compute the eigenvalues of the staggered microscale model~\cref{eqs:stgrdFDEs_gnrcWve} over the full-domain, we follow well established  analysis (\cites[e.g.,][pp.138--139]{Hinch2020_ThnkBfreYuCmpte_APrldeToCmpttnlFldDynmcs}[]{Griffiths2011_TrvlngWvAnlyssPDENmrclAnlytclMthdsMtlbMple})---which we summarise here to compare its results with the patch scheme.
We substitute into the staggered finite difference system~\cref{eqs:stgrdFDEs_gnrcWve} an arbitrary Fourier mode\footnote{Throughout this article, \(i\)~denotes the micro-grid index in \(x\)-direction, whereas \(\i := \sqrt{-1}\) is the imaginary unit.} \(h, u, v \propto  \exp[ \i (k_x {\cij i} \delta  + k_y {\cij j} \delta )  + \lambda t]\) with real wavenumber~$(k_x, k_y)$ and complex growth rate~\(\lambda\).  Then factoring the complex exponentials leads to the eigensystem
\begin{equation} \label{eqn:eigSys_m1_gLnrWve}
  \mathbf{J}_\mu
  \begin{bmatrix}
    h \\
    u \\
    v
  \end{bmatrix}
  =
  \lambda_\mu 
  \begin{bmatrix}
    h \\
    u \\
    v
  \end{bmatrix},
\end{equation}
where the subscript~\(\mu\) denotes the full-domain microscale staggered model (to distinguish it from the patch scheme), and where the Jacobian  
\begin{align*}
&  
  \mathbf{J}_\mu := 
  \begin{bmatrix}
  0 
    & - \i \sin{\left (k_{x} \delta \right )} / \delta
    & - \i \sin{\left (k_{y} \delta \right )} / \delta
    \\
  - \i \sin{\left (k_{x} \delta \right )} / \delta
    & - c_D - c_V \omega_{\mu0}^2
    & 0
    \\
  - \i \sin{\left ( k_{y} \delta \right )} / \delta
    & 0
    & - c_D - c_V \omega_{\mu0}^2
\end{bmatrix},
\end{align*}
with frequency \(\omega_{\mu0} := \sqrt{ \sin^2{\left(k_x \delta \right) }/\delta^2 + \sin^2{\left(k_y \delta \right) }/\delta^2 } \) corresponding to ideal waves (\(c_D = c_V = 0\)), for the {staggered grid} microscale model.
For each wavenumber~\((k_x,k_y)\), the three eigenvalues of the \(3 \times 3\) Jacobian~\(\mathbf{J}_\mu\) are (one real and a complex conjugate pair),
\begin{align}& \label{eqn:eig_muA}
\lambda_\mu  =
\begin{cases}
  - c_D + c_V \omega_{\mu0}^2\,,\\
  - \left( c_D + c_V \omega_{\mu0}^2 \right)/2
    \pm \i \sqrt{
      \omega_{\mu0}^2
      - \left[ \left( c_D + c_V \omega_{\mu0}^2 \right)/2 \right]^2
      } \,.
\end{cases}
\end{align}
We assess the accuracy of a staggered patch scheme by comparing the eigenvalues of its macroscale modes (corresponding to small wavenumbers) with these eigenvalues~\eqref{eqn:eig_muA} for the corresponding small wavenumber modes in the full-domain microscale model.

Second, we compute the eigenvalues of the staggered patch scheme~\cref{eqs:FDEs_pN_gnrcWve}.
The dynamics of a staggered patch scheme reflects the coupled dynamics at two different length scales.
One length scale is due to the microscale interactions within the patches, and the other length scale is due to the macroscale patch coupling.
Hence, for eigenvalue analysis of a staggered patch scheme, the general mode needs to include spatial structures of the two scales: microscale within a patch, and macroscale across the patches.
\begin{itemize}
    \item The staggered patch scheme~\cref{eqs:FDEs_pN_gnrcWve} is invariant to translations in space by multiples of the macroscale spacing~\(2\Delta\).
    Hence the macroscale spatial structure can be expressed by complex exponential factors \(\exp[ \i (k_x {\cIJ I} \Delta + k_y {\cIJ J} \Delta)]\). The indices~\({\cIJ I, J}\) increment by multiples of two due to the translation symmetry.
    \item But the microscale spatial structure within the patches is not translationally invariant on the microscale: for example, nodes near a patch edge are different to nodes on the patch interior.
    Hence, in contrast to the full-domain microscale Fourier mode~\(h, u, v \propto  \exp[ \i (k_x {\cij i} \delta  + k_y {\cij j} \delta )  + \lambda t]\), with constant amplitudes for~\(h,u,v\), a staggered patch scheme mode must allow for a general microscale spatial structure within the patches.
    We represent such sub-patch microscale structure by the variation with the sub-patch micro-grid indices~\(i,j\) of the fields~\(h_{{\cij i,j}}^{{\cpq p,q}}\), \(u_{{\cij i,j}}^{{\cpq p,q}}\), \(v_{{\cij i,j}}^{{\cpq p,q}}\) for each patch within a macro-cell ({\cpq orange} squares in \cref{fig:indexing_patchGrid}); here
    \({{\cpq p,q}} \in \{0,1\}\) is the local sub-macro-cell patch index.
\end{itemize}
Thus, for the eigenvalue analysis of a staggered patch scheme, we consider an arbitrary staggered patch scheme Fourier mode with real \emph{macroscale wavenumber}~$(k_x, k_y)$, over an infinite staggered grid of patches (\(N \to \infty\) in \cref{fig:PtchGrd_n2t0}),
\begin{subequations} \label{eqs:FrrMde_ptch_gLnrWve}
\begin{align}
  \hNode\,\hVarP(t)
    &= h_{{\cij i,j}}^{{\cpq p,q}}(t) 
      \exp[ \i (k_x {\cIJ I} \Delta + k_y {\cIJ J} \Delta)] ,\\
  \uNode\,\uVarP(t)
    &= u_{{\cij i,j}}^{{\cpq p,q}}(t)  
      \exp[ \i (k_x {\cIJ I} \Delta + k_y {\cIJ J} \Delta)] ,\\
  \vNode\,\vVarP(t)
    &= \,v_{{\cij i,j}}^{{\cpq p,q}}(t)  
      \exp[ \i (k_x {\cIJ I} \Delta + k_y {\cIJ J} \Delta)] ,
\end{align}
\end{subequations}
for the sub-patch micro-grid indices~\(({\cij i}, {\cij j})\), global macroscale patch indices~\(({\cIJ I}, {\cIJ J})\), and local sub-macro-cell patch index~\(({\cpq p,q})\) in the patch scheme~\cref{eqs:FDEs_pN_gnrcWve}.
The local sub-cell patch indices, \({\cpq p} := {\cIJ I} \bmod 2\) and \({\cpq q} := {\cIJ J} \bmod 2\).
\Cref{fig:indexConvention} illustrates these indices for an example case.
In a finite space domain the global macroscale patch index~\({\cIJ I,J} \in \{0,1, \ldots, N-1\}\). 
For algebraic eigenvalue analysis, we invoke an infinite domain, \(N\to\infty\)\,, for which the macroscale patch index~\({\cIJ I, J} \in\ZZ\), the integers.

A staggered patch grid over the \(2\pi \times 2\pi\) domain (\(L = 2\pi\)) with inter-patch distance~\(\Delta\),  resolves macroscale wavenumbers ranging from the smallest resolved wavenumber~\(k_x = k_y = 2\pi/L = 1\) to the largest resolved wavenumber~\(\lfloor 2\pi / (4 \Delta) \rfloor = \lfloor \pi / (2 \Delta) \rfloor\).
Macroscale modes corresponding to these relatively small macroscale wavenumbers~\(|k_x|,|k_y| \leq \pi/(2\Delta)\) are all linearly independent on the macroscale patch grid, but not so for higher wavenumbers.
Linear combinations of these linearly independent modes give general solutions to the patch systems~\eqref{eqs:FDEs_pN_gnrcWve}, 
so we restrict the modes to this wavenumber range.

In the patch scheme mode~\cref{eqs:FrrMde_ptch_gLnrWve}, the time-dependent microscale spatial structure~\(h_{{\cij i,j}}^{{\cpq p,q}}(t)\), \(u_{{\cij i,j}}^{{\cpq p,q}}(t)\), \(v_{{\cij i,j}}^{{\cpq p,q}}(t)\) is \emph{spatially modulated} by the macroscale wave form \(\exp[\i (k_x {\cIJ I} \Delta + k_y {\cIJ J} \Delta)]\).
The microscale structure~\(h_{{\cij i,j}}^{{\cpq p,q}}(t)\), \(u_{{\cij i,j}}^{{\cpq p,q}}(t)\), \(v_{{\cij i,j}}^{{\cpq p,q}}(t)\) depends on the sub-macro-cell patch index~\({{\cpq p,q}}\) and the sub-patch micro-grid node index~\({\cij i,j}\), but not on the global patch index~\({\cIJ I,J}\).

In eigenvalue analysis of a patch scheme, the time-dependent microscale structure~\(h_{{\cij i,j}}^{{\cpq p,q}}(t)\), \(u_{{\cij i,j}}^{{\cpq p,q}}(t)\), \(v_{{\cij i,j}}^{{\cpq p,q}}(t)\) in the patch scheme Fourier mode~\cref{eqs:FrrMde_ptch_gLnrWve} corresponds to the interior nodes of all patches in any one macro-cell, indicated by the orange square in \cref{fig:PtchGrd_n2t0}.
Collecting the interior values of all patches in the macro-cell into a vector gives the \emph{state vector}~\(\vec{x}^i\); the superscript~\({(\argdot)^i}\) is not an index or exponent, instead, a qualifier denoting the patch interior nodes of one macro-cell (we use~\({(\argdot)^I}\) to denote the interior nodes of all the macro-cells in the full system~\cref{eqn:dynSys_pN_gnrcWve}).
For the example patch grid~\#79985 of \cref{fig:PtchGrd_n1t0_huvx,fig:mltiLyrEdgeNdsForPatchGrid}, the total number of patch interior nodes per macro-cell, that is the size of~\(\vec{x}^i\),
\begin{equation} \label{eqn:exprsn_ni_gLnrWve}
  n_p^i = 9n^2/4 - 4n + 2,
\end{equation}
where \(n\) is the number of sub-patch micro-grid intervals.
For example, for the cases of \(n = 6,\, 10,\, 14\) sub-patch micro-grid intervals, \(n_p^i = 59,\, 187,\, 387\) respectively.
For the specific staggered patch grid~\#79985 with~\(n=6\) sub-patch micro-grid intervals, the state vector containing \(59\)~elements is
\begin{equation} \label{eqs:stateVect_p1_n6} 
  \vec{x}^i = \parbox[t]{0.8\linewidth}{\mathcode`\,="213B$
  \big(
  h^{{\cpq0,0}}_{{\cij1,1}},
  h^{{\cpq0,0}}_{{\cij1,3}},
  h^{{\cpq0,0}}_{{\cij1,5}},
  h^{{\cpq0,0}}_{{\cij3,1}},
  h^{{\cpq0,0}}_{{\cij3,3}},
  h^{{\cpq0,0}}_{{\cij3,5}},
  h^{{\cpq0,0}}_{{\cij5,1}},
  h^{{\cpq0,0}}_{{\cij5,3}},
  h^{{\cpq0,0}}_{{\cij5,5}},
  u^{{\cpq0,0}}_{{\cij2,1}},
  u^{{\cpq0,0}}_{{\cij2,3}},
  u^{{\cpq0,0}}_{{\cij2,5}},
  u^{{\cpq0,0}}_{{\cij4,1}},
  u^{{\cpq0,0}}_{{\cij4,3}},
  u^{{\cpq0,0}}_{{\cij4,5}},
  v^{{\cpq0,0}}_{{\cij1,2}},
  v^{{\cpq0,0}}_{{\cij1,4}},
  v^{{\cpq0,0}}_{{\cij3,2}},
  v^{{\cpq0,0}}_{{\cij3,4}},
  v^{{\cpq0,0}}_{{\cij5,2}},
  v^{{\cpq0,0}}_{{\cij5,4}},
  h^{{\cpq0,1}}_{{\cij1,2}},
  h^{{\cpq0,1}}_{{\cij1,4}},
  h^{{\cpq0,1}}_{{\cij3,2}},
  h^{{\cpq0,1}}_{{\cij3,4}},
  h^{{\cpq0,1}}_{{\cij5,2}},
  h^{{\cpq0,1}}_{{\cij5,4}},
  u^{{\cpq0,1}}_{{\cij2,2}},
  u^{{\cpq0,1}}_{{\cij2,4}},
  u^{{\cpq0,1}}_{{\cij4,2}},
  u^{{\cpq0,1}}_{{\cij4,4}},
  v^{{\cpq0,1}}_{{\cij1,1}},
  v^{{\cpq0,1}}_{{\cij1,3}},
  v^{{\cpq0,1}}_{{\cij1,5}},
  v^{{\cpq0,1}}_{{\cij3,1}},
  v^{{\cpq0,1}}_{{\cij3,3}},
  v^{{\cpq0,1}}_{{\cij3,5}},
  v^{{\cpq0,1}}_{{\cij5,1}},
  v^{{\cpq0,1}}_{{\cij5,3}},
  v^{{\cpq0,1}}_{{\cij5,5}},
  h^{{\cpq1,0}}_{{\cij2,1}},
  h^{{\cpq1,0}}_{{\cij2,3}},
  h^{{\cpq1,0}}_{{\cij2,5}},
  h^{{\cpq1,0}}_{{\cij4,1}},
  h^{{\cpq1,0}}_{{\cij4,3}},
  h^{{\cpq1,0}}_{{\cij4,5}},
  u^{{\cpq1,0}}_{{\cij1,1}},
  u^{{\cpq1,0}}_{{\cij1,3}},
  u^{{\cpq1,0}}_{{\cij1,5}},
  u^{{\cpq1,0}}_{{\cij3,1}},
  u^{{\cpq1,0}}_{{\cij3,3}},
  u^{{\cpq1,0}}_{{\cij3,5}},
  u^{{\cpq1,0}}_{{\cij5,1}},
  u^{{\cpq1,0}}_{{\cij5,3}},
  u^{{\cpq1,0}}_{{\cij5,5}}, 
  v^{{\cpq1,0}}_{{\cij2,2}},
  v^{{\cpq1,0}}_{{\cij2,4}},
  v^{{\cpq1,0}}_{{\cij4,2}},
  v^{{\cpq1,0}}_{{\cij4,4}}
  \big).
  $}
\end{equation}

Recall that the patch coupling (spectral throughout this article) gives the edge values of each patch as a function of the aggregate values of all the patches.
Hence a patch coupling gives edge values of all the patches in a macro-cell in terms of the substituted macroscale mode wavenumber.
Collecting the edge values of all patches in the macro-cell into a vector gives the \emph{edge vector}~\(\vec{x}^e\); the superscript~\({(\argdot)^e}\) is not an index or exponent, instead, a qualifier denoting the edge nodes of one macro-cell (we use~\(E\) to denote the edge nodes of all the macro-cells in the full system~\cref{eqn:dynSys_pN_gnrcWve}).
For the particular staggered patch grid~\#79985 in \cref{fig:PtchGrd_n2t0}, the total number of patch edge nodes per macro-cell for weakly damped linear wave, that is the size of the edge vector~\(\vec{x}^e\), is
  \(n^e_p = 18n - 16\)\,,
where \(n\) is the number of sub-patch micro-grid intervals: for \(n = 6,\, 10,\, 14\), the vector size \(n^e_p = 92,\, 164,\, 236\), respectively.
This~\(n_e^i\) is only for the \(48\) centred patch grids with three patches per macro-cell (\cref{fig:PtchGrd_n2t0}), the \(72\) centred patch grids with four patches per macro-cell have more edge nodes.

Substituting the macroscale Fourier mode~\cref{eqs:FrrMde_ptch_gLnrWve}, and the coupled patch edge values computed by the patch coupling, into the staggered patch scheme~\cref{eqs:FDEs_pN_gnrcWve} of the wave \textsc{pde}s~\cref{eqs:PDEs_gLnrWve} gives the time evolution of a staggered patch scheme as a linear dynamical system
\begin{equation} \label{eqn:dynSys_p1_gLnrWve}
  \dby{\vec{x}^i}{t} = \vec{F}(\vec{x}^i; \vec{x}^e(\vec{x}^i))
  = \mathbf{J}_p  \vec{x}^i ,
\end{equation}
where the Jacobian \(\mathbf{J}_p  := \doh{\vec{F}}{\vec{x}^i}\)\,.
This system is parametrised by the macroscale wavenumber~\((k_x, k_y)\).
There are \(n_p^i\)~variables corresponding to each macroscale wavenumber~\((k_x, k_y)\).
The state vector~\(\vec{x}^I\) of the full-size staggered patch scheme dynamical system~\cref{eqn:dynSys_pN_gnrcWve} contain interior values of all the patches in a staggered patch grid, but the state vector~\(\vec{x}^i\) of the staggered patch scheme dynamical system~\cref{eqn:dynSys_p1_gLnrWve} for one macroscale wavenumber~\((k_x, k_y)\) contain interior values of only one macro-cell.
Hence, we call equation~\cref{eqn:dynSys_p1_gLnrWve} the \emph{one-cell system} of a staggered patch scheme.
Similarly, we call the Jacobian~\(\mathbf{J}_p\) as the \emph{one-cell Jacobian} of a staggered patch scheme.

The \(n_p^i \times n_p^i\) one-cell Jacobian~\(\mathbf{J}_p \) depends on the physical parameters~\(c_D\) and~\(c_V\), discretisation parameters~\(n\), \(\delta\), and~\(\Delta\), and the macroscale wavenumber~\((k_x,k_y)\).
For the example patch grid~\#79985 in \cref{fig:PtchGrd_n2t0}, for \(n=6\) sub-patch micro-grid intervals, the one-cell Jacobian~\(\mathbf{J}_p \) is a \(59 \times 59\) sparse matrix with at most only~\(318\) nonzero elements of its \(3481\)~elements.
Such sparsity holds irrespective of the particular patch coupling interpolation.
%

By considering all macroscale wavenumbers~\((k_x,k_y)\), the one-cell Jacobian~\(\mathbf{J}_p \) provides a complete general solution for every initial condition applied to a given patch scheme.
Let the eigenvalues and eigenvectors of~\(\mathbf{J}_p \) be denoted by~\(\lambda_l\) and~\(\ev_l\) respectively.
Then by the linearity of~\eqref{eqn:dynSys_p1_gLnrWve}, a general solution is a linear combination of the set of~\(\ev_l\exp(\lambda_lt)\).
So stability is assured if the eigenvalues~\(\{\lambda_l\}\) all have nonpositive real parts.
The few eigenvectors~\(\ev_l\) that have spatial structures that are slowly varying within the patches are the macroscale modes in the patch scheme.
It is the eigenvalues of these modes that we compare for accuracy with the small wavenumber eigenvalues of the full-domain system.
The other eigenvectors~\(\ev_l\), those with significant sub-patch structure, are micro-scale modes modulated over a macroscale.
They are not physically significant, and so we do not compare their eigenvalues with that of the underlying full-domain system, although we do require that they be stable.
It is in this way that the eigenvalues of the Jacobian~\(\mathbf{J}_p \) characterise the stability and accuracy of a given patch scheme.

We attempted to derive algebraic expressions for the eigenvalues of a one-cell Jacobian~\(\mathbf{J}_p \), through various algebraic simplification strategies, in the three Computer Algebra Systems (\textsc{cas}) SymPy, Reduce, and Maple.
For general macroscale wavenumber~\((k_x,k_y)\), the three \textsc{cas} packages failed to compute (no results in 48~hours) the analytic eigenvalues of the \(59 \times 59\) Jacobian (\(n=6\)).
That is, unlike the expression~\cref{eqn:eig_muA} for the full-domain micro-scale models, expressing the eigenvalues algebraically appears infeasible for the patch scheme.
Hence, we compute the patch scheme eigenvalues~\(\lambda_p \) by numerically evaluating the one-cell Jacobian for a range of values of the macroscale wavenumber~\((k_x,k_y)\), physical parameters~\(c_D\) and~\(c_V\), and discretisation parameters~\(n\), \(\delta\), and~\(\Delta\).
%
%

\Cref{ssc:onlyCntrdPtchGrdsAreAcrte,ssc:onlyPtchGrdsWthAlSymtrcPtchsAreStble} assess the staggered patch grids by comparing eigenvalues~\(\lambda_p \) of a patch scheme one-cell Jacobian~\(\mathbf{J}_p\) with the eigenvalues~\(\lambda_\mu\) of the full domain microscale model, for all the macroscale wavenumbers~\((k_x,k_y)\) resolved by a finite macroscale staggered grid.

\subsection{Only patch grids with all symmetric patches are stable}
\label{ssc:onlyPtchGrdsWthAlSymtrcPtchsAreStble}

\Cref{ssc:extndTheStgrdGrdToMltsclePtchSchme} introduces compatible multiscale staggered patch grids for wave-like systems and shows that there are~\(167\,040\) such compatible 2D patch grids.
Among these staggered patch grids, \(120\) of them are \emph{centred staggered patch grids} (e.g., the two patch grids in \cref{fig:PtchGrd_n1t0_huvx_huvh}) whose sub-patch micro-grids all have a centre node and include a centred-patch corresponding to each of the state variables.
That is, the centred patch grids consist of each of the {\cij green} highlighted patches in \cref{fig:allMicroGrids} and the fourth patch may also be one such centred patch or an empty patch.
The remaining~\(166\,920\) are \emph{non-centred staggered patch grids}.
This subsection studies the dependence of the patch scheme stability on the patch grid geometry, the \emph{geometry-stability study}, for all these~\(167\,040\) compatible 2D patch grids for the weakly damped \text{linear wave \pde{}s~\cref{eqs:PDEs_gLnrWve}.}

We define a \emph{staggered patch scheme to be stable} when the maximum real part of its \emph{numerically computed} eigenvalues~\(\max\Re(\lambda_p)< 10^{-5}\) (to allow for tiny round-off errors).
And herein we define a \emph{staggered patch grid to be stable} if the spectral staggered patch scheme is stable over that grid.

\begin{figure}
  \centering
  \caption{\label{fig:maxRe_ptchGrids_idlWve_Spectral_N10_n6_r0p1}%
    Maximum real parts of the spectral patch scheme eigenvalues~\(\max\Re(\lambda_p)\) over all~\(83\,520\) compatible 2D staggered patch grids with \(N=10\), \(n=6\), \(r=0.1\) for ideal wave \pde{}s~\cref{eqs:PDEs_gLnrWve} with~\(c_D = c_V = 0\).
    Among the non-centred patch grids (red dots), \(82\,896\)~are unstable (\(\max\Re(\lambda_p) > 10^{-5}\)), \(564\)~are stable (\(\max\Re(\lambda_p) < 10^{-5}\)).
    All~\(60\) centred staggered patch grids (blue circles) are stable.
    }
  \includegraphics[scale=1]{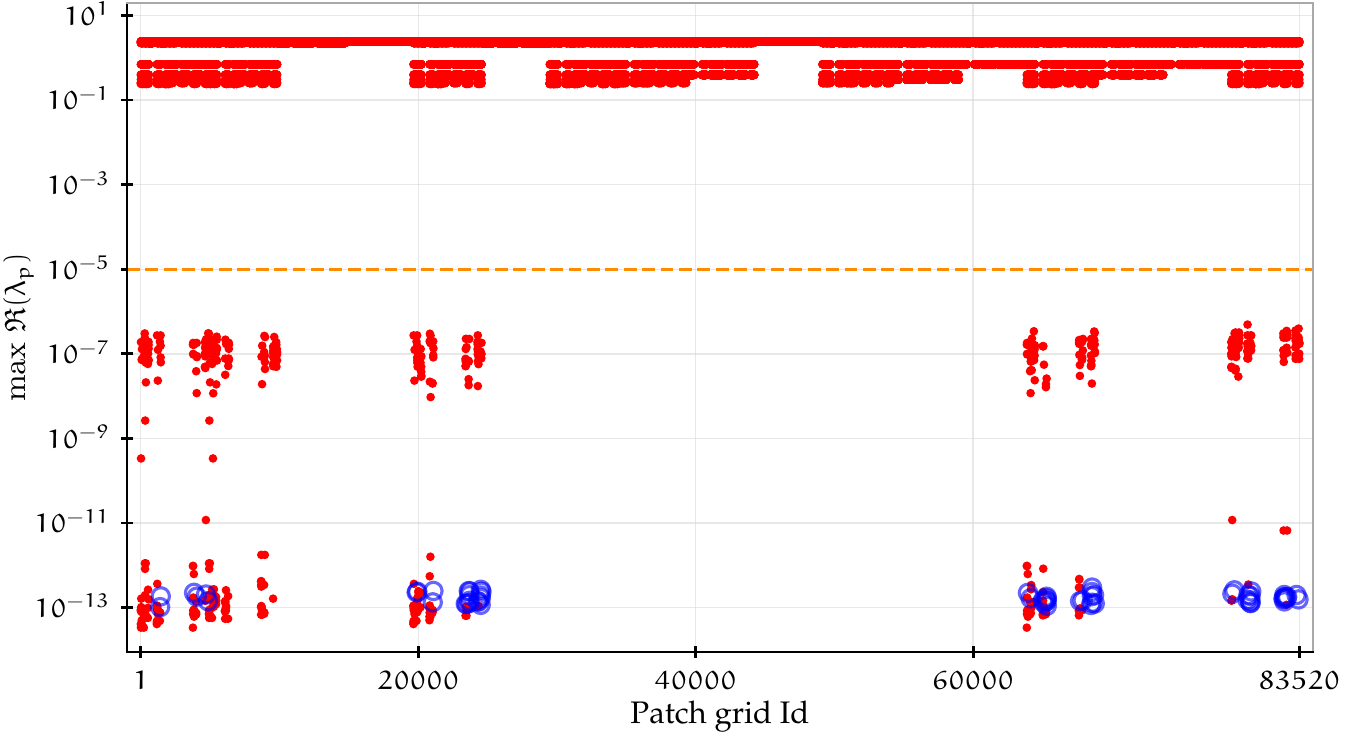}
\end{figure}

We analyse the stability of all the~\(167\,040\) staggered patch grids separately for two groups based on the number of sub-patch micro-grid intervals~\(n\).
First, we analyse the \(83\,520\) patch grids with odd~\(n/2\), taking~\(n=6\) as a representative case; second the \(83\,520\) patch grids with odd~\(n/2\), taking~\(n=4\) as a representative case.
\Cref{fig:maxRe_ptchGrids_idlWve_Spectral_N10_n6_r0p1} plots the maximum real parts of the spectral patch scheme eigenvalues~\(\max\Re(\lambda_p)\) over all~\(83\,520\) compatible 2D staggered patch grids with \(N=10\), \(n=6\), \(r=0.1\) for ideal wave \pde{}s~\cref{eqs:PDEs_gLnrWve} with~\(c_D = c_V = 0\).
The red dots are for non-centred patch grids, and the blue circles are for centred patch grids.
In each patch grid, each of the four patches in a macro-cell is chosen from any of the \(16\) micro-grids in \cref{fig:allMicroGrids_n6} or an empty patch, excluding the all-empty patch.
Hence, \cref{fig:maxRe_ptchGrids_idlWve_Spectral_N10_n6_r0p1} is representative of the staggered patch grids with odd~\(n/2\).
In \cref{fig:maxRe_ptchGrids_idlWve_Spectral_N10_n6_r0p1}, among the non-centred patch grids (red dots), only \(564\) of them are stable according to \(\max\Re(\lambda_p) < 10^{-5}\), all the remaining \(82\,896\) of the non-centred patch grids are unstable with \(\max\Re(\lambda_p) > 10^{-5}\).
That is, \(99.32\%\) of the non-centred patch grids are unstable.
On the other hand all~\(60\) of the centred staggered patch grids (blue circles) are \text{stable with~\(\max\Re(\lambda_p) < 3 \cdot 10^{-13}\).}
That is, \emph{almost all non-centred patch grids lead to unstable patch schemes but all centred patch grids constitute stable patch schemes}.

In \cref{fig:maxRe_ptchGrids_idlWve_Spectral_N10_n6_r0p1}, all~\(564\) stable non-centred patch grids (red dots with \(\max\Re(\lambda_p) < 10^{-5}\)), contain only symmetric patches~\(uuvv\), \(hhvv\), \(uuhh\), and \(hhhh\) from \cref{fig:allMicroGrids_n6} or an empty patch.
All~\(60\) stable centred patch grids (blue circles with \(\max\Re(\lambda_p) < 3 \cdot 10^{-13}\)), contain only \(h\)-, \(u\)-, and \(v\)-centred patches from \cref{fig:allMicroGrids_n6} or an empty patch, which are also symmetric patches.
Thus, \cref{fig:maxRe_ptchGrids_idlWve_Spectral_N10_n6_r0p1} shows that among all~\(83\,520\) staggered patch grids, \emph{only patch grids containing symmetric patches are stable}.
On the other hand, only \(5^4 - 1 = 624\) staggered patch grids are possible using only the symmetric patches~\(uuvv\), \(hhvv\), \(uuhh\), and \(hhhh\) from \cref{fig:allMicroGrids_n6} or an empty patch, excluding the all-empty patch.
Among these~\(60\) are centred patch grids and the remaining~\(564\) are non-centred patch grids.
\Cref{fig:maxRe_ptchGrids_idlWve_Spectral_N10_n6_r0p1} shows all these~\(624\) patch grids with symmetric patches are stable.
This shows that \emph{all patch grids containing only symmetric patches are stable}.
That is, \emph{using only symmetric patches is a necessary and sufficient condition for patch scheme stability}.

As a representative of the staggered patch grids with even~\(n/2\), since even cases are potentially different to the odd cases, we plotted the maximum real parts of the eigenvalues for~\(n=4\) (not included).
%
%
The plot was qualitatively the same as \cref{fig:maxRe_ptchGrids_idlWve_Spectral_N10_n6_r0p1} except that the clusters of stable patch grids occur around different patch grid Ids than in \cref{fig:maxRe_ptchGrids_idlWve_Spectral_N10_n6_r0p1}.
Such a difference is because the centredness depends upon whether~\(n/2\) is odd or even ({\cij green} highlighted names in \cref{fig:allMicroGrids_n4,fig:allMicroGrids_n6}).
For \(n=4\), the spectral patch is stable, \(\max\Re(\lambda_p) < 10^{-5}\), for all the \(624\) staggered patch grids that contain only symmetric patches~\(uuvv\), \(hhvv\), \(uuhh\), and \(hhhh\) from \cref{fig:allMicroGrids_n4} or an empty patch.
That plot of \(\max\Re(\lambda_p)\) for~\(n=4\) also confirms the same two conclusions in the preceding two paragraphs, including the same number of stable and unstable patch grids.
Thus, for both the cases of~\(n/2\) being odd or even, among all the possible~\(167\,040\) compatible 2D staggered patch grids, only \(624 + 624 = 1248\) patch grids (\(0.75\%\)) are stable. 
%

Including the dissipation in the linear wave \pde{}s~\cref{eqs:PDEs_gLnrWve} by nonzero dissipation coefficients~\(c_D, c_V\) stabilize the patch schemes by pushing the positive real parts of the eigenvalues to negative values.
Hence all these~\(1248\) patch grids are also stable for weakly damped linear waves.

\begin{figure}
\caption{\label{fig:ptchGrid_cmplxEigs_unstable1}%
  Non-centred staggered patch grid~\#55420 over which spectral patch scheme is unstable for ideal wave \pde{}s~\cref{eqs:PDEs_gLnrWve} with~\(c_D = c_V = 0\).
}
\centering
\begin{subfigure}[t]{0.475\textwidth}
  \centering
  \caption{\label{fig:ptchGrd_unstable1}%
    %
    This non-centred patch grid~\#55420 has the non-centred patch~\(uhvh\) from \cref{fig:allMicroGrids_n6} in place of the \(h\)-patch of the patch grid~\#79985 (\cref{fig:PtchGrd_n1t0_huvx}); that is the only difference.
  }%
  \includegraphics[scale=1]{figs/PatchGrid_idlWve_N6n6_n1t0_55420_unstable1_\nodeShape}
\end{subfigure}%
\quad
\begin{subfigure}[t]{0.475\textwidth}
  \centering
  \caption{\label{fig:cmplxEigs_unstable1}%
    Eigenvalues of spectral patch scheme over patch grid~\#55420 has large positive real parts, \(\max\Re(\lambda_p) = 0.3\).
  }%
  \includegraphics[scale=1]{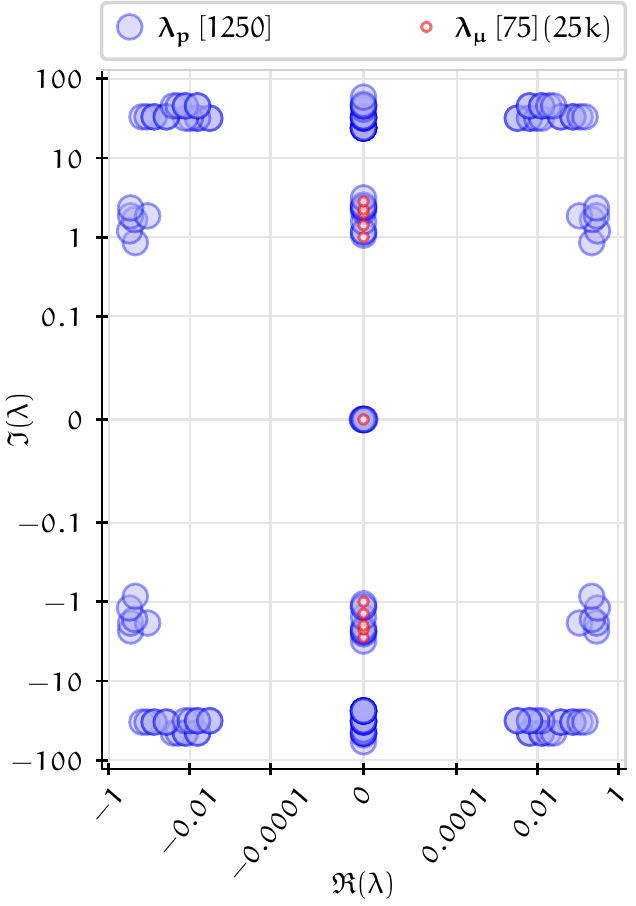}
\end{subfigure}
\end{figure}

\begin{figure}
\caption{\label{fig:ptchGrid_cmplxEigs_unstable3}%
  Non-centred patch grid~\#56236 leads to more unstable modes compared to patch grid~\#55420 of \cref{fig:ptchGrd_unstable1} for ideal waves.
}
\centering
\begin{subfigure}[t]{0.475\textwidth}
  \centering
  \caption{\label{fig:ptchGrd_unstable3}%
    Replacing \(h\)-, \(u\)-, \(v\)-patches of the centred patch grid~\#79985 (\cref{fig:PtchGrd_n1t0_huvx}) with non-centred patches~\(uhvh\), \(huvh\), \(uhhv\) (\cref{fig:allMicroGrids_n6}) respectively, gives this non-centred patch grid~\#56236.
  }%
  \includegraphics[scale=1]{figs/PatchGrid_idlWve_N6n6_n1t0_56236_unstable3_\nodeShape}
\end{subfigure}%
\quad
\begin{subfigure}[t]{0.475\textwidth}
  \centering
  \caption{\label{fig:cmplxEigs_unstable3}%
    Eigenvalues of spectral patch scheme over patch grid~\#56236 has large positive real parts, \(\max\Re(\lambda_p) = 0.3\).
  }%
  \includegraphics[scale=1]{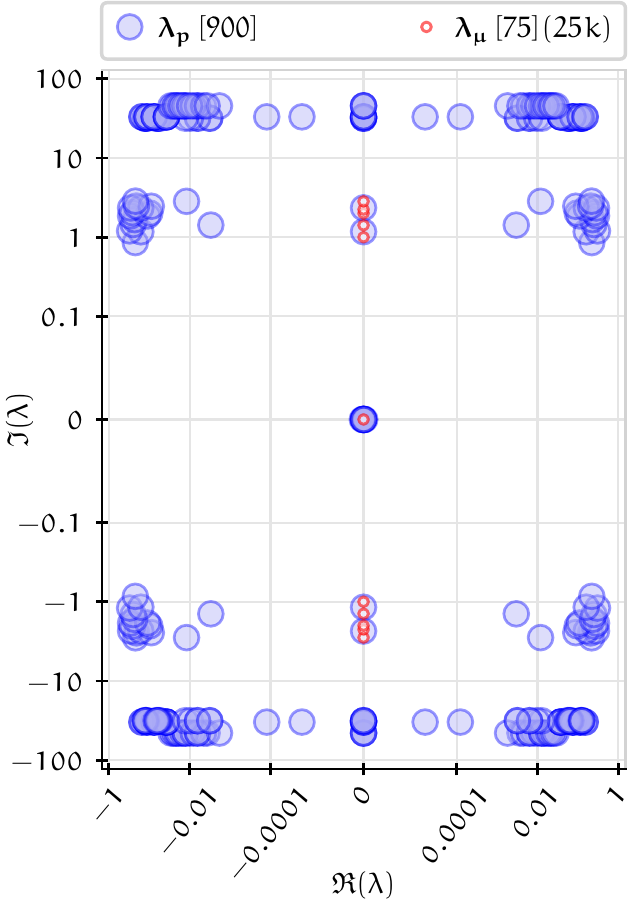}
\end{subfigure}
\end{figure}

This paragraph discusses three examples of unstable non-centred staggered patch grids and their eigenvalues.
Firstly, \cref{fig:ptchGrid_cmplxEigs_unstable1} shows the non-centred patch grid~\#55420 and the corresponding eigenvalues of the spectral patch scheme for discretisations of wave \pde{}s~\cref{eqs:PDEs_gLnrWve} with~\(c_D = c_V = 0\).
%
The patch grid~\#55420 (\cref{fig:ptchGrd_unstable1}) has the non-centred patch~\(uhvh\) (\cref{fig:allMicroGrids_n6}) in place of the \(h\)-patch of~\#79985 (\cref{fig:PtchGrd_n1t0_huvx}).
That one change to a non-centred patch makes the patch grid~\#55420 unstable with large real part eigenvalues~\(\max\Re(\lambda_p) = 0.3\) (\cref{fig:cmplxEigs_unstable1}).
There are~\(80\) unstable modes with \(\max\Re(\lambda_p) > 10^{-5}\) in \cref{fig:cmplxEigs_unstable1} for the non-centred patch grid~\#55420.
Secondly,  replacing the \(u\)-patch with the non-centred patch~\(huvh\) (\cref{fig:allMicroGrids_n6}) in the patch grid~\#55420 of \cref{fig:ptchGrd_unstable1}, gives the non-centred patch grid~\#56287 which we found has \(176\)~unstable modes (plot not included).
Lastly, replacing also the \(v\)-patch in~\#56287  with the non-centred patch~\(uhhv\) (\cref{fig:allMicroGrids_n6}) gives the patch grid~\#56236 in \cref{fig:ptchGrd_unstable3} whose patches are all non-centred.
\cref{fig:cmplxEigs_unstable3} plots the eigenvalues for this non-centred patch grid~\#56236: it shows \(272\)~unstable modes with \(\max\Re(\lambda_p) > 10^{-5}\).
It appears that {the more non-centred patches in a staggered patch grid, the more modes are unstable}.

\subsection{Only centred patch grids are accurate}
\label{ssc:onlyCntrdPtchGrdsAreAcrte}

This subsection shows that all the \(60\)~centred staggered patch grids give accurate patch schemes, but none of the stable non-centred patch grids in \cref{fig:maxRe_ptchGrids_idlWve_Spectral_N10_n6_r0p1} is accurate.
We first discuss discretisations of the ideal wave \pde{}s~\cref{eqs:PDEs_gLnrWve} (\(c_D = c_V = 0\)), and second explore weakly damped linear waves with nonzero~\(c_D, c_V\).
This subsection also illustrates several features and correlations in the structure of the eigenvalues of the patch scheme for weakly damped linear waves.


\begin{figure}
  \centering
    \begin{minipage}[t][][c]{0.5\linewidth}
      \captionof{table}{%
        Slow macroscale eigenvalues~\(\lambda_p\) of the patch scheme~\cref{eqs:FDEs_pN_gnrcWve}, \(\left|\Im(\lambda_p)\right| < 10\) in \cref{fig:cmplxEigs_idlWve_Spectral_N10_n6_r0p01_huvx_oce}, agree exactly (within~\(10^{-12}\)) with the corresponding eigenvalues~\(\lambda_\mu\) of the full-domain microscale model~\eqref{eqs:stgrdFDEs_gnrcWve}.
        All the distinct ideal wave frequencies~\(\omega_{\mu 0}\) (red circles in \cref{fig:cmplxEigs_idlWve_Spectral_N10_n6_r0p01_huvx_oce}) correspond to the resolved macroscale wavenumbers~\(k_x, k_y = 0, \pm 1, \pm 2\) on the patch grid~\#79985 with~\(N = 10\).
        Error~\(\epsilon_k = \max_i(|\lambda_{p,i} - \lambda_\mu|\) for each wavenumber.
      }
      \label{tbl:macroEigsAgreeExactly}
    \end{minipage}
    \quad
    \begin{minipage}[t][][c]{0.45\linewidth}
      \begin{tabular}{rlc}
        \hline
        \(\omega_{\mu 0}\) &
          \((k_x, k_y)\) &
          \(\epsilon_k\)
          \\
          \hline
          \(-2.8284\) & \((-2, -2)\) & \(2.0\cdot 10^{-13}\)  \\
          \(-2.2361\) & \((-2, -1)\) & \(2.6\cdot 10^{-13}\)  \\
          \(-2.0000\) & \((-2, 0)\) & \(2.0\cdot 10^{-13}\)  \\
          \(-1.4142\) & \((-1, -1)\) & \(1.6\cdot 10^{-13}\)  \\
          \(-1.0000\) & \((-1, 0)\) & \(2.7\cdot 10^{-13}\)  \\
          \(0.0000\) & \((0, 0)\) & \(2.9\cdot 10^{-14}\)  \\
          \(1.0000\) & \((1, 0)\) & \(2.3\cdot 10^{-13}\)  \\
          \(1.4142\) & \((1, 1)\) & \(1.2\cdot 10^{-13}\)  \\
          \(2.0000\) & \((2, 0)\) & \(1.3\cdot 10^{-13}\)  \\
          \(2.2361\) & \((2, 1)\) & \(2.1\cdot 10^{-13}\)  \\
          \(2.8284\) & \((2, 2)\) & \(1.9\cdot 10^{-13}\)  \\
        \hline
        \end{tabular}
    \end{minipage}
  \end{figure}    

This paragraph illustrates the accuracy of the centred staggered patch grid~\#79985 for ideal waves.
%
For non-dissipative ideal wave \pde{}s~\cref{eqs:PDEs_gLnrWve} (\(c_D = c_V = 0\)), \cref{fig:cmplxEigs_idlWve_Spectral_N10_n6_r0p01_huvx_oce} compares eigenvalues~\(\lambda_p\) of the patch scheme~\cref{eqs:FDEs_pN_gnrcWve} with eigenvalues~\(\lambda_p\) of the full-domain microscale model~\eqref{eqs:stgrdFDEs_gnrcWve} with the same micro-grid interval~\(\delta = 2\pi/3000\).
\Cref{fig:cmplxEigs_idlWve_Spectral_N10_n6_r0p01_huvx_oce} shows that the slow macroscale eigenvalues, \(\left|\Im(\lambda_p)\right| < 10\), of the patch scheme agree to graphical accuracy with the eigenvalues~\(\lambda_\mu\) of the full-domain microscale staggered model.
\Cref{tbl:macroEigsAgreeExactly} shows this agreement is exact (within~\(10^{-12}\)).
First column in \cref{tbl:macroEigsAgreeExactly} lists all the distinct ideal wave frequencies \(\omega_{\mu0} = \Im(\lambda_\mu)\) corresponding to the red circles in \cref{fig:cmplxEigs_idlWve_Spectral_N10_n6_r0p01_huvx_oce}, using the eigenvalue expression~\cref{eqn:eig_muA}.
The second column lists all the resolved macroscale wavenumbers~\(k_x, k_y = 0, \pm 1, \pm 2\) on the patch grid~\#79985 with~\(N = 10\) that corresponds to these~\(\omega_{\mu0}\).
For each of these wavenumbers, there are several repeated patch scheme eigenvalues~\(\lambda_p\) (blue circles in \cref{fig:cmplxEigs_idlWve_Spectral_N10_n6_r0p01_huvx_oce}).
For each of these wavenumbers, the third column lists the errors~\(\epsilon_k = \max_i(|\lambda_{p,i} - \lambda_\mu|)\) between the eigenvalues of the patch scheme and that of the full-domain model.
The small overall maximum error of~\(\epsilon_k = 2.7\cdot 10^{-13}\) (within numerical roundoff error) shows that slow macroscale eigenvalues~\(\lambda_p\) and \(\lambda_\mu\) agree exactly. 
That is, the centred staggered patch grid~\#79985 gives accurate multiscale spectral patch schemes for ideal waves.

\begin{figure}
  \centering
  \caption{\label{fig:imEig_ptchGrds_sym_idlWve_Spectral_N10_n6_r0p1}%
    Imaginary parts of spectral patch scheme macroscale eigenvalues (\(\left|\Im(\lambda_p)\right| < 10\)) for all \(624\)~stable staggered patch grids that contain only symmetric patches, with \(N=10\), \(n=6\), \(r=0.1\) for ideal wave \pde{}s~\cref{eqs:PDEs_gLnrWve} with~\(c_D = c_V = 0\).
    For all \(60\)~centred patch grids \(\Im(\lambda_p)\) (blue circles) agree with \(\Im(\lambda_\mu)\) (black lines) of full-domain microscale model~\eqref{eqs:stgrdFDEs_gnrcWve}.
    For none of the \(564\)~non-centred patch grids \(\Im(\lambda_p)\) (red dots) agree with all the \(\Im(\lambda_\mu)\).
  }
  \includegraphics[scale=1]{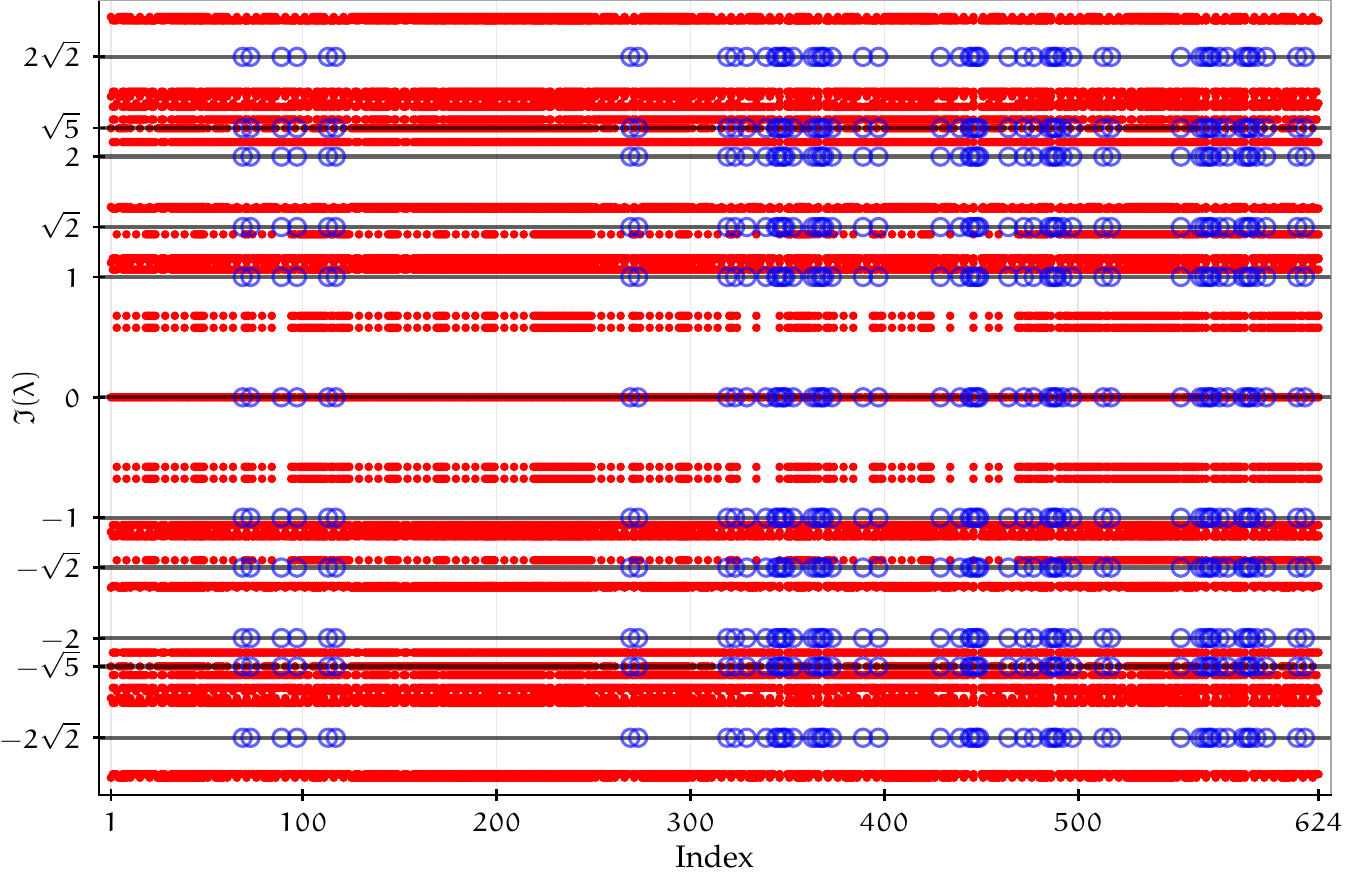}
\end{figure}

This paragraph shows that all the centred staggered patch grids give accurate multiscale spectral patch schemes for ideal waves.
The previous \Cref{ssc:onlyPtchGrdsWthAlSymtrcPtchsAreStble} shows that all  \(624\)~patch grids containing only symmetric patches are stable (i.e., in \cref{fig:allMicroGrids_n6}, patches \(uuvv\), \(hhvv\), \(uuhh\), \(hhhh\), or empty).
The negligible~\(\max\Re(\lambda_p)\) of these stable schems agrees with the zero real parts in the eigenvalue expression~\cref{eqn:eig_muA} for full-domain microscale model of ideal waves (\(c_D = c_V = 0\)).
That all eigenvalues~\(\lambda_p\) are pure imaginary (to numerical error) shows that the staggered patch grid preserves waves without introducing any significant artificial dissipation.
Hence, to assess the accuracy of all those~\(624\) stable staggered patch grids, it is sufficient to compare the imaginary parts of the eigenvalues.
\Cref{fig:imEig_ptchGrds_sym_idlWve_Spectral_N10_n6_r0p1} plots the imaginary parts of the slow macroscale eigenvalues, those for which \(\left|\Im(\lambda_p)\right| < 10\), for all the~\(624\) stable staggered patch grids that contain only symmetric patches.
The horizontal black lines in \cref{fig:imEig_ptchGrds_sym_idlWve_Spectral_N10_n6_r0p1}  correspond to the ideal wave frequencies \(\omega_{\mu0} = \Im(\lambda_\mu) = \sqrt{ \sin^2{\left(k_x \delta \right) }/\delta^2 + \sin^2{\left(k_y \delta \right) }/\delta^2} \approx (0,\, \pm 1,\, \pm \sqrt{2},\, 2,\, \pm\sqrt{5},\, \pm 2\sqrt{2})\) from the eigenvalue expression~\cref{eqn:eig_muA} for macroscale wavenumbers~\(k_x, k_y = 0, \pm 1, \pm 2\) and \(\delta = 2\pi/300\).  
For all \(60\)~centred patch grids (blue circles)  \(\Im(\lambda_p)\) agree exactly with \(\Im(\lambda_\mu)\) (black lines) of full-domain microscale model~\eqref{eqs:stgrdFDEs_gnrcWve} (within numerical roundoff errors).
On the other hand, for none of the \(564\)~non-centred patch grids \(\Im(\lambda_p)\) (red dots) agree with all the \(\Im(\lambda_\mu)\) (black lines).
That is, \emph{for ideal waves, all the \(60\)~centred staggered patch grids give accurate patch schemes, but none of the stable non-centred patch grids is accurate}.

Several features in the structure of patch scheme eigenvalues become clearer when plotted for cases of the weakly damped linear wave \pde{}s~\cref{eqs:PDEs_gLnrWve} with nonzero dissipation coefficients~\(c_D, c_V\).
The rest of this subsection discusses the structure of patch scheme eigenvalues and establishes the patch scheme accuracy for weakly damped linear wave \pde{}s~\cref{eqs:PDEs_gLnrWve}.
Unless otherwise noted, we show and discuss the case of small bed drag with \(c_D = 10^{-6}\) and small viscous dissipation with \(c_V = 0.0001\).

\begin{figure}
\centering
\caption{\label{fig:eigPlt_gLnrWve_Spectral_N14_n6_r0p1_cD1e-06_cV0p0001_allK}%
    Eigenvalues for weakly damped linear waves~(\(c_D = 10^{-6}\), \(c_V = 0.0001\)) over a \(2\pi \times 2\pi\) periodic domain with the micro-grid interval \(\delta = 2\pi/420\).
    For macroscale modes with \(\Re(\lambda) > -0.01\) and \(\left|\Im(\lambda)\right| < 10\), the eigenvalues~\(\lambda_p\) of the spectral staggered patch scheme (\(N=14\), \(n=6\), \(r=0.1\)) agree exactly with eigenvalues~\(\lambda_\mu\) of the staggered full-domain microscale model.
    }
  \includegraphics[scale=1]{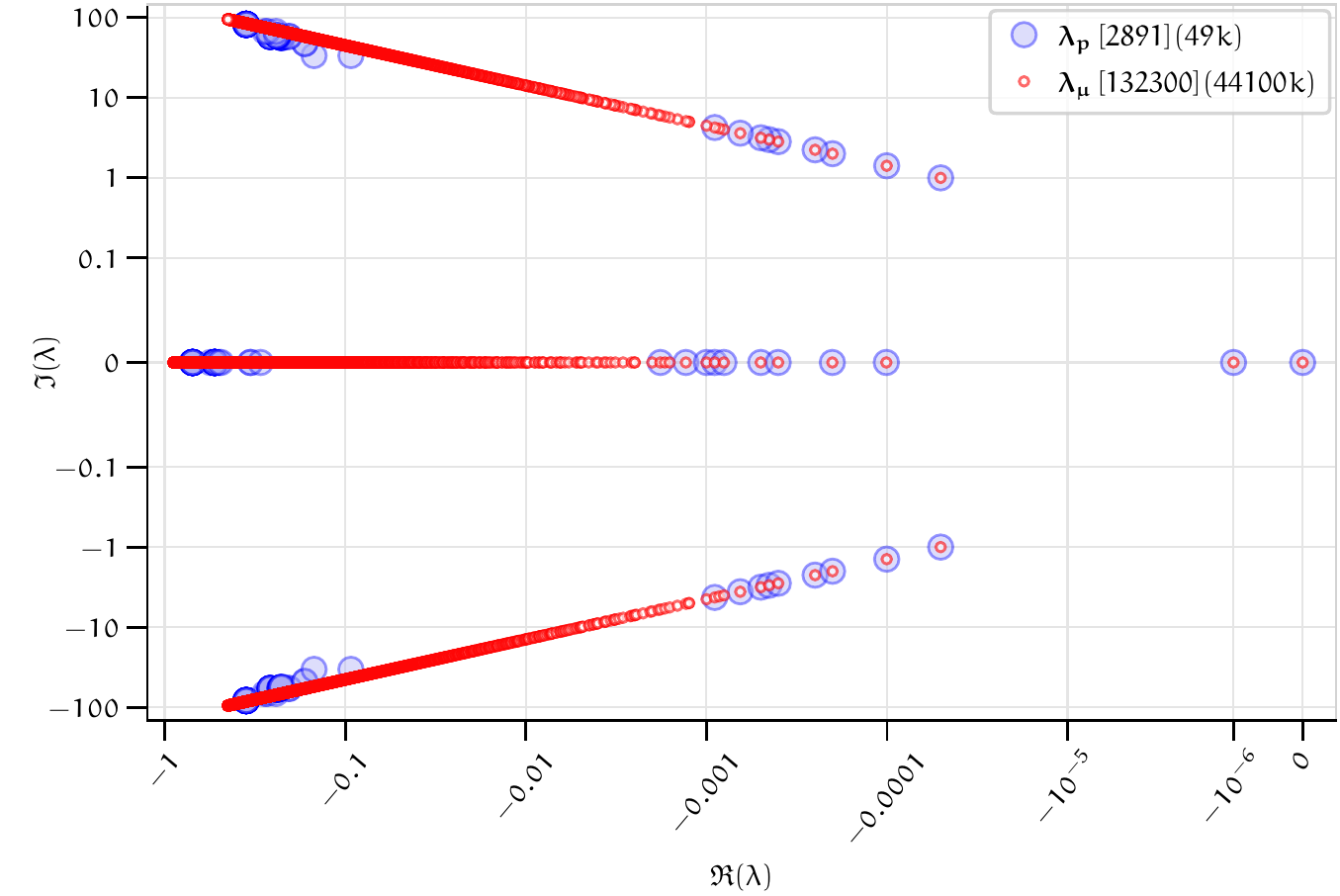}
\end{figure}

\Cref{fig:eigPlt_gLnrWve_Spectral_N14_n6_r0p1_cD1e-06_cV0p0001_allK} compares the patch scheme eigenvalues~\(\lambda_p \) with the eigenvalues~\(\lambda_\mu\) of the full-domain microscale model~\cref{eqs:stgrdFDEs_gnrcWve} for the weakly damped linear wave.
The eigenvalues~\(\lambda_p\) in \cref{fig:eigPlt_gLnrWve_Spectral_N14_n6_r0p1_cD1e-06_cV0p0001_allK} are for the spectral patch scheme over the centred staggered patch grid~\#79985 (\cref{fig:PtchGrd_n1t0_huvx}) over a non-dimensional \(2\pi \times 2 \pi\) domain (\(L = 2\pi\)) with \(14 \times 14\) (\(N = 14\)) macro-grid intervals, \(6 \times 6\) (\(n = 6\)) sub-patch micro-grid intervals, and patch scale ratio~\(r = 0.1\).
Hence, the macro-grid interval~\(\Delta = L/N = 2\pi/14\) ({\cIJ violet} macro-grid) and the sub-patch micro-grid interval~\(\delta = l / n = 2 L r / (Nn) = 2\pi/420 \approx 0.015\) ({\cij green} sub-patch micro-grid).
A staggered patch grid with \(N=14\)~macro-grid intervals resolves \(N^2/4 = 49\)~macroscale wavenumbers~\((k_x, k_y)\) which is same as the number of macro-cells.
As per the count~\cref{eqn:exprsn_ni_gLnrWve}, such a staggered patch grid has~\(n_p^i = 9n^2/4 - 4n + 2 = 59\) patch interior nodes per macro-cell forming the one-cell state vector~\(\vec{x}^i\) in the expression~\cref{eqs:stateVect_p1_n6}.
Hence, for each macroscale wavenumber~\((k_x, k_y)\) resolved by the staggered patch grid, the \(n_p^i \times n_p^i\) one-cell Jacobian~\(\mathbf{J}_p \) in the linear system~\cref{eqn:dynSys_p1_gLnrWve} gives~\(n_p^i = 59\) eigenvalues~\(\lambda_p \).
Numerically evaluating~\(\mathbf{J}_p \) for all \(49\)~macroscale wavenumbers~\((k_x, k_y)\) then gives \(49 \cdot 59 = 2891\) patch scheme eigenvalues~\(\lambda_p \).
As required, this matches the total number of patch interior nodes~\(n_p^I\) (the number of state variables) from the expression~\cref{eqn:exprsn_nI_gLnrWve} for this case of \(N = 14\) and \(n = 6\).

The staggered full-domain microscale model for \cref{fig:eigPlt_gLnrWve_Spectral_N14_n6_r0p1_cD1e-06_cV0p0001_allK} is over the full-domain micro-grid (\cref{fig:flDmnGrd_clctd_stgrd}, right)  whose micro-grid interval is of the same size as the sub-patch micro-grid interval \(\delta = 2\pi/420\) of the patch scheme for the figure.
Such a full-domain micro-grid has \(n = L/\delta = 420\) micro-grid intervals and resolves \(n^2/4 = 44\,100\) wavenumbers~\((k_x, k_y)\).
%
%
For each wavenumber~\((k_x, k_y)\) resolved on the staggered full-domain micro-grid, the \(3 \times 3\) Jacobian~\(\mathbf{J}_\mu\) in the linear system~\cref{eqn:eigSys_m1_gLnrWve} gives \(3\)~eigenvalues~\(\lambda_\mu\), as in the expression~\cref{eqn:eig_muA}.
Thus, evaluating the expression~\cref{eqn:eig_muA} for all the \(44\,100\)~wavenumbers gives \(44\,100 \cdot 3 = 132\,300\) full-domain microscale eigenvalues~\(\lambda_\mu\).
As required, this matches the total number of full-domain microscale nodes~\(3n^2/4\) (the number of state variables) for \(n = 420\) (\cref{sec:cllctdAndStgrdGrdsFrWvelkeSystms}).

In \Cref{fig:eigPlt_gLnrWve_Spectral_N14_n6_r0p1_cD1e-06_cV0p0001_allK}, the eigenvalues~\(\lambda_\mu\) of the staggered full-domain microscale model (small red circles) denote the modes corresponding to all \(44\,100\)~resolved wavenumbers, \emph{ranging uniformly across the spatial scales}, from the smallest resolved wavenumber~\(k_x = k_y = 2\pi/L = 1\) to the largest resolved wavenumber~\(2\pi / (4 \delta) = 105\).
On the other hand, the eigenvalues~\(\lambda_p \) of the multiscale staggered patch scheme (large blue circles) correspond to two qualitatively different modes at two different spatial scales, pure \emph{macroscale modes} and \emph{sub-patch microscale modes}, with a gap between their spatial scales.
This gap in the spatial scales also leads to a gap in their time scales, evident in the \emph{spectral gap} where there are no patch scheme eigenvalues roughly between \(-0.1 \lesssim \Re(\lambda_p)  \lesssim -0.001\) in \cref{fig:eigPlt_gLnrWve_Spectral_N14_n6_r0p1_cD1e-06_cV0p0001_allK}.
\begin{itemize}
  \item \emph{Macroscale modes} are those patch scheme modes (eigenvectors) that have macroscale structure (spatial variation) with little microscale structure.
  The eigenvalues corresponding to these macroscale modes are \emph{macroscale eigenvalues} (\(\Re(\lambda_p)\gtrsim-0.001\) in \cref{fig:eigPlt_gLnrWve_Spectral_N14_n6_r0p1_cD1e-06_cV0p0001_allK}).
  \item \emph{Microscale modes} are those patch schemes modes that have significant microscale structure irrespective of whether it is modulated by some macroscale structure.
  The eigenvalues corresponding to these microscale modes are \emph{microscale eigenvalues} (\(\Re(\lambda_p)\lesssim-0.1\) in \cref{fig:eigPlt_gLnrWve_Spectral_N14_n6_r0p1_cD1e-06_cV0p0001_allK}).
\end{itemize}
Plots (omitted in this article) of corresponding eigenvectors confirm the above division into macroscale and microscale modes \parencite[\S3.2.6]{Divahar2022_AcrteMltiscleSmltnOfWveLkeSystms} .
The macroscale eigenvalues~\(\lambda_p \) on the right five clusters correspond to all the \(N^2/4 = 49\)~macroscale wavenumbers resolved by the staggered patch grid, ranging from the smallest macroscale wavenumber~\(k_x = k_y = 2\pi/L = 1\) to the largest resolved wavenumber~\(\lfloor 2\pi / (4 \Delta) \rfloor = \lfloor N/4 \rfloor = 3\); that is the macroscale wavenumber~\(k_x, k_y \in \{ -3, -2, -1, 0, 1, 2, 3 \}\).
The microscale eigenvalues~\(\lambda_p \) on the left three clusters correspond to sub-patch microscale modes modulated over the macroscale, whose modulated wavenumbers range from the smallest sub-patch microscale wavenumber~\(k_x = k_y = 2\pi/l = 2\pi/(2 r L/N) = 70\) to the largest sub-patch microscale wavenumber~\(2\pi / (4 \delta) = 105\).
Hence, the gap between the largest resolved macroscale wavenumber and the smallest sub-patch microscale wavenumber spanning \((\lfloor N/4 \rfloor, 2\pi/(2 r L/N)) = (3, 70)\), is the gap in spatial scales that a patch scheme does not resolve.
In \cref{fig:eigPlt_gLnrWve_Spectral_N14_n6_r0p1_cD1e-06_cV0p0001_allK}, the \emph{spectral gap} where there are no patch scheme eigenvalues roughly between \(-0.1 \lesssim \Re(\lambda_p ) \lesssim -0.001\), corresponds to these spatial scales unresolved by the patch scheme.

\begin{figure}
\centering
  \caption{\label{fig:eigPlt_gLnrWve_Spectral_N14_n6_r0p1_cD1e-06_cV0p0001}%
    \emph{All} macroscale eigenvalues~\(\lambda_p\) agree exactly with eigenvalues~\(\lambda_\mu\).
    Same as \cref{fig:eigPlt_gLnrWve_Spectral_N14_n6_r0p1_cD1e-06_cV0p0001_allK} except here the eigenvalues~\(\lambda_\mu\) of staggered full-domain microscale model are plotted only for those \(N^2/4 = 49\)~macroscale wavenumbers resolved by the staggered patch grid.
    }
  \includegraphics[scale=1]{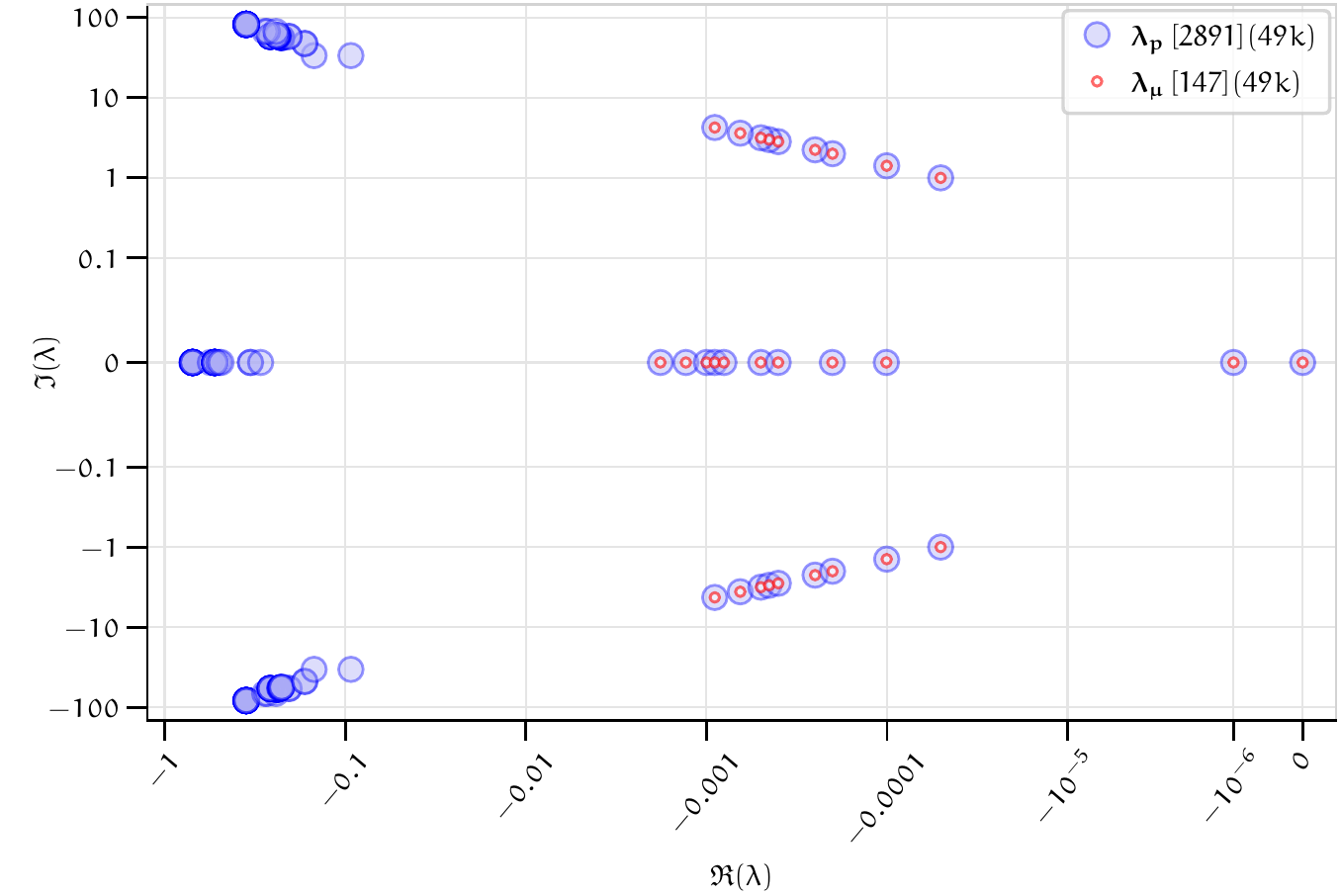}
\end{figure}

\cref{fig:eigPlt_gLnrWve_Spectral_N14_n6_r0p1_cD1e-06_cV0p0001} plots the eigenvalues~\(\lambda_\mu\) corresponding to only those \(N^2/4\)~macroscale wavenumbers resolved by the staggered patch grid.
%
Similar to \cref{fig:cmplxEigs_idlWve_Spectral_N10_n6_r0p01_huvx_oce} for ideal wave, \cref{fig:eigPlt_gLnrWve_Spectral_N14_n6_r0p1_cD1e-06_cV0p0001} for weakly damped linear wave clearly shows that \emph{all} macroscale eigenvalues~\(\lambda_p \) with  \(\Re(\lambda_p) > -0.01\) and \(\left|\Im(\lambda_p) \right| < 10\) agree exactly (within numerical roundoff errors) with eigenvalues~\(\lambda_\mu\) of the full-domain microscale model for the same macroscale wavenumbers.
That is, the equation-free multiscale spectral \emph{patch scheme on the staggered patch grid~\#79985 (\cref{fig:PtchGrd_n2t0}) accurately resolves the macroscale weakly damped linear waves}.

As \cref{sec:cllctdAndStgrdGrdsFrWvelkeSystms} discusses, it is the microscale computational model that a multiscale patch scheme seeks to predict accurately, not directly the \textsc{pde} model.
That is, how well the full-domain microscale models~\cref{eqs:clctdFDEs_gnrcWve,eqs:stgrdFDEs_gnrcWve} predict the solutions of the \textsc{pde}s~\cref{eqs:PDEs_gnrcWve} is a peripheral issue.

\begin{figure}
  \centering
  \caption{\label{fig:reEig_ptchGrds_sym_gLnrWve_Spectral_N10_n6_r0p1_cD1e-06_cV0p0001}%
    Real parts of spectral patch scheme macroscale eigenvalues (\(\Re(\lambda_p) > -0.01\) and \(\left|\Im(\lambda_p)\right| < 10\)) for all \(624\)~stable staggered patch grids that contain only symmetric patches, with \(N=10\), \(n=6\), \(r=0.1\) for weakly damped linear wave \pde{}s~\cref{eqs:PDEs_gLnrWve} with~\(c_D = 10^{-6}\), \(c_V = 0.0001\).
    For all \(60\)~centred patch grids \(\Re(\lambda_p)\) (blue circles) agree with \(\Re(\lambda_\mu)\) (black lines) of full-domain microscale model~\eqref{eqs:stgrdFDEs_gnrcWve}.
    For none of the \(564\)~non-centred patch grids \(\Re(\lambda_p)\) (red dots) agree with all the \(\Re(\lambda_\mu)\).
  }
  \includegraphics[scale=1]{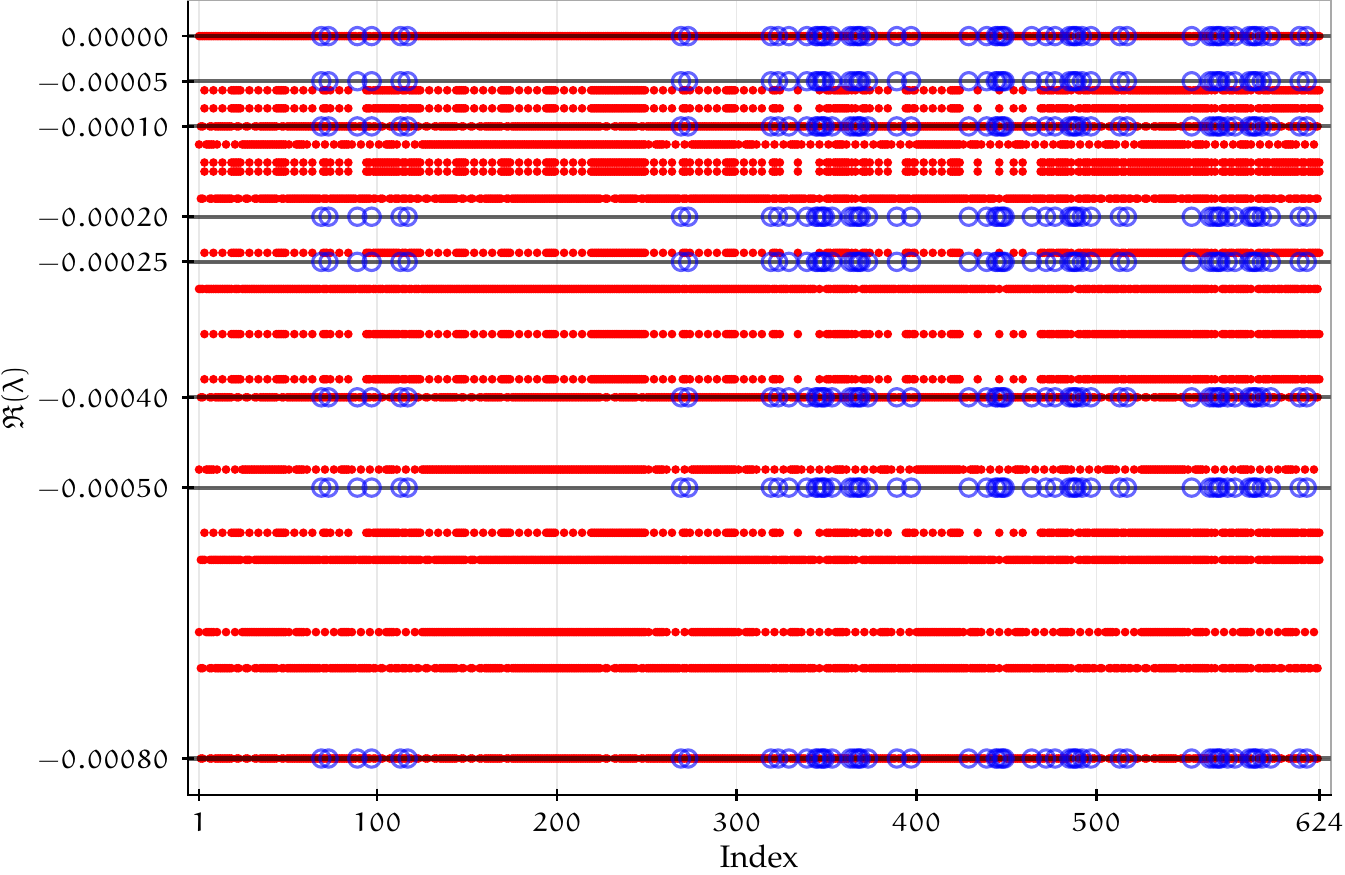}
\end{figure}

\emph{Centred staggered patch grids give accurate multiscale spectral patch scheme for weakly damped linear waves.}
For weakly damped linear wave \pde{}s~\cref{eqs:PDEs_gLnrWve} with~\(c_D = 10^{-6}\), \(c_V = 0.0001\), \cref{fig:eigPlt_gLnrWve_Spectral_N14_n6_r0p1_cD1e-06_cV0p0001} illustrate that the macroscale eigenvalues (\(\Re(\lambda_p) > -0.01\) and \(\left|\Im(\lambda_p)\right| < 10\)) of the patch scheme\#79985 agree exactly with the eigenvalues~\(\lambda_\mu\) of the staggered full-domain microscale model~\eqref{eqs:stgrdFDEs_gnrcWve}.
We now explore all of the \(624\)~stable staggered patch grids, and the imaginary and real parts of their eigenvalues.
First, a plot of imaginary parts of the macroscale eigenvalues for all the~\(624\) stable staggered patch grids for weakly damped linear waves (not included herein) is visually nearly identical to \cref{fig:imEig_ptchGrds_sym_idlWve_Spectral_N10_n6_r0p1} for the ideal wave.
Quantitatively, for all \(60\)~centred patch grids, \(\Im(\lambda_p)\)~agrees to numerical error with~\(\Im(\lambda_\mu)\) of the full-domain microscale model.
But none of the \(564\)~non-centred patch grids~\(\Im(\lambda_p)\) agree with all the~\(\Im(\lambda_\mu)\).
Secondly, \cref{fig:reEig_ptchGrds_sym_gLnrWve_Spectral_N10_n6_r0p1_cD1e-06_cV0p0001} plots the real parts of the macroscale eigenvalues (red dots and blue circles) for all the~\(624\) stable staggered patch grids (\cref{fig:maxRe_ptchGrids_idlWve_Spectral_N10_n6_r0p1}) that contain only symmetric patches~\(uuvv\), \(hhvv\), \(uuhh\), and \(hhhh\) from \cref{fig:allMicroGrids_n6} or an empty patch.
In \cref{fig:reEig_ptchGrds_sym_gLnrWve_Spectral_N10_n6_r0p1_cD1e-06_cV0p0001}, the horizontal black lines correspond to \(\Re(\lambda_\mu)\) from the eigenvalues~\cref{eqn:eig_muA} for macroscale wavenumbers~\(k_x, k_y = 0, \pm 1, \pm 2\) with \(\delta = 2\pi/300\).
\Cref{fig:reEig_ptchGrds_sym_gLnrWve_Spectral_N10_n6_r0p1_cD1e-06_cV0p0001} illustrates that for all \(60\)~centred patch grids, \(\Re(\lambda_p)\)~(blue circles) agree exactly with \(\Re(\lambda_\mu)\) (black lines) of full-domain microscale model~\eqref{eqs:stgrdFDEs_gnrcWve} (within numerical roundoff errors).
On the other hand, for none of the \(564\)~non-centred patch grids do all~\(\Re(\lambda_p)\) (red dots) agree with the \(\Re(\lambda_\mu)\) (black lines).
That is, just as for ideal waves, \emph{for weakly damped linear waves, all the \(60\)~centred staggered patch grids give accurate patch schemes, but none of the stable non-centred patch grids are accurate}.

All the accuracy studies reported in the preceding paragraphs of this subsection are for \(n=6\) sub-patch micro-grid intervals, which are representative of the staggered patch grids with odd~\(n/2\).
A similar study for~\(n=4\) (not included herein), as a representative of patch grids with even~\(n/2\), also confirms that among the~\(83\,520\) possible 2D staggered patch grids only the \(60\)~centred staggered patch grids give accurate patch schemes.
Thus, among \emph{all these possible~\(167\,040\) compatible 2D staggered patch grids, only \(120\) centred patch grids (\(0.07\%\)) are both stable and accurate}.

\section{Conclusion}

Accurate numerical simulation of wave-like systems is challenging, especially over long simulation times.
Such challenges become even more intricate for large-scale simulations reliant on modelling small-scale features, especially in multiple dimensions.
For wave-like systems in full-domain modelling, a common strategy for accurate and robust numerical schemes is to use staggered grids (\cref{sec:cllctdAndStgrdGrdsFrWvelkeSystms}).
%
%
For high accuracy and to preserve much of the wave characteristics, \emph{this article extends the concept of staggered grids in full-domain modelling to multidimensional multiscale modelling} (\cref{sec:stgrdPtchGrdsFrEqtnFreMltscleMdlngOfWvs}).
We developed and analysed all the possible~\(167\,040\) 2D staggered patch grids that are \emph{geometrically compatible} for wave-like systems (\cref{ssc:extndTheStgrdGrdToMltsclePtchSchme}).
For first-order \pde{}s, the patch schemes over these staggered patch grids interpolate field values to the patch edges.
However, higher-order spatial derivatives (e.g., for diffusion) and some nonlinear terms require interpolation of additional layers of edge values (\cref{ssc:mltiLyrEdgeNdsFrHghrOrdrSptlDrvtvs}).
Extending to higher dimensions should be straightforward.

Almost all such multiscale staggered patch grids lead to unstable and/or inaccurate equation-free multiscale patch schemes.
We identified~\(120\) staggered patch grids that constitute stable and accurate multiscale schemes.
For representative physical and discretisation parameters, via eigenvalues and wave frequencies, we demonstrate the stability, accuracy, and wave-preserving characteristic of the centred multiscale staggered grids for weakly damped linear waves (\cref{sec:mltscleStgrdPtchGrdFrWklyDmpdLnrWvs}).
The geometry-stability study (\cref{ssc:onlyPtchGrdsWthAlSymtrcPtchsAreStble}) shows that among all the compatible~\(167\,040\) compatible patch grids, only \(1248\)~patch grids (\(0.75\%\)) whose patches are all symmetric are stable.
None of the non-centred patch grids that are stable (\cref{ssc:onlyPtchGrdsWthAlSymtrcPtchsAreStble}) is accurate (\cref{ssc:onlyCntrdPtchGrdsAreAcrte}).
Thus, among all the possible~\(167\,040\) compatible 2D staggered patch grids, only~\(120\) centred patch grids (\(0.07\%\)) are both stable and accurate.

%
This article develops \(120\) centred multiscale staggered grids and demonstrates their stability, accuracy, and wave-preserving characteristic for multiscale modelling of weakly damped linear waves for representative physical and discretisation parameters.
%
%
But most characteristics of the developed multiscale staggered grids must also hold in general for multiscale modelling of many complex spatio-temporal physical phenomena such as the general computational fluid dynamics.

Our next article on multiscale simulation of large-scale linear waves explores two families of patch schemes of \cref{ssc:ptchCplngCnctsTheScls} in greater detail over a wider range of parameters for their stability, accuracy, consistency, and their practical insensitivity to numerical roundoff errors.
Subsequent articles will explore staggered patch schemes for nonlinear wave \pde{}s, specifically for viscous thin fluid flows and for turbulent shallow water flows.

\end{document}